\documentclass[a4paper,11pt]{article}
\usepackage{latexsym}
\usepackage{amsmath}
\usepackage{amsthm}

\usepackage{graphicx}
\usepackage{txfonts}
\usepackage{bm}
\usepackage{color}	
\usepackage{epic,eepic}

\newcommand{\RED}[1]{{\color{red}#1}} 
\newcommand{\BLU}[1]{{\color{blue}#1}} 

  \renewcommand{\RED}[1]{{\color{black}#1}} 
  \renewcommand{\BLU}[1]{{\color{black}#1}}

\usepackage{geometry}
\numberwithin{equation}{section}

\newtheorem{theorem}{Theorem}[section]
\newtheorem{proposition}[theorem]{Proposition}
\newtheorem{corollary}[theorem]{Corollary}
\newtheorem{lemma}[theorem]{Lemma}
\newtheorem{claim}[theorem]{Claim}

\newtheorem{remark}{Remark}[section]
\newtheorem{example}{Example}[section]

\newcommand{\memo}[1]{{\bf \small \RED{[MEMO:}} \BLU{#1} \ {\bf \small \RED{:end]} }}  
   \renewcommand{\memo}[1]{}           
\newcommand{\OMIT}[1]{{\bf [OMIT:} #1 \ {\bf --- end OMIT] }}  
   \renewcommand{\OMIT}[1]{}            

\newcommand{\RR}{{\bf R}}
\newcommand{\ZZ}{{\bf Z}}

\newcommand{\finbox}{\hspace*{\fill}$\rule{0.17cm}{0.17cm}$}
\newcommand{\finboxHere}{\ $\rule{0.17cm}{0.17cm}$}

\newcommand{\odotZ}{\overset{....}}

\newcommand{\llceil}{\bigg\lceil} 
\newcommand{\rrceil}{\bigg\rceil} 

\newcommand{\Proof}{\noindent {\bf Proof.  }}

\usepackage{lineno}
\newcommand*\patchAmsMathEnvironmentForLineno[1]{
  \expandafter\let\csname old#1\expandafter\endcsname\csname #1\endcsname
  \expandafter\let\csname oldend#1\expandafter\endcsname\csname end#1\endcsname
  \renewenvironment{#1}
     {\linenomath\csname old#1\endcsname}
     {\csname oldend#1\endcsname\endlinenomath}}
\newcommand*\patchBothAmsMathEnvironmentsForLineno[1]{
  \patchAmsMathEnvironmentForLineno{#1}
  \patchAmsMathEnvironmentForLineno{#1*}}
\AtBeginDocument{
\patchBothAmsMathEnvironmentsForLineno{equation}
\patchBothAmsMathEnvironmentsForLineno{align}
\patchBothAmsMathEnvironmentsForLineno{flalign}
\patchBothAmsMathEnvironmentsForLineno{alignat}
\patchBothAmsMathEnvironmentsForLineno{gather}
\patchBothAmsMathEnvironmentsForLineno{multline}
}

\begin{document}

\title{
Fair Integral Network Flows%
\footnote{
This is a revised version of 
``Discrete Decreasing Minimization, Part III:  Network Flows,'' 2019.
http://arxiv.org/abs/1907.02673
}
}

\author{Andr\'as Frank%
\thanks{MTA-ELTE Egerv\'ary Research Group,
Department of Operations Research, E\"otv\"os University, P\'azm\'any
P.~s.~1/c, Budapest, Hungary, H-1117. 
e-mail:  {\tt frank@cs.elte.hu}. 
ORCID: 0000-0001-6161-4848.
The research was partially supported by the
National Research, Development and Innovation Fund of Hungary
(FK\_18) -- No. NKFI-128673.
}
 \ \ and \ 
{Kazuo Murota%
\thanks{
The Institute of Statistical Mathematics,
Tokyo 190-8562, Japan; 
Faculty of Economics and Business Administration,
Tokyo Metropolitan University,
Tokyo 192-0397, Japan,
e-mail:  {\tt murota@tmu.ac.jp}. 
ORCID: 0000-0003-1518-9152.
The research was supported by 
JSPS KAKENHI Grant Numbers JP26280004, JP20K11697. 
 }}}


\date{September 2020 / January 2022 / April 2022}

\maketitle

\begin{abstract}
A strongly polynomial algorithm is developed 
for finding an integer-valued feasible $st$-flow of given flow-amount
which is decreasingly minimal 
on a specified subset $F$ of edges 
in the sense that the largest flow-value on $F$ is as small as possible,
within this, the second largest flow-value on $F$ is as small as possible, 
within this, the third largest flow-value on $F$ is as small as possible, 
and so on.
A characterization of the set of these
$st$-flows gives rise to an algorithm to compute a cheapest
$F$-decreasingly minimal integer-valued feasible $st$-flow of given flow-amount.  
Decreasing minimality is a possible formal way to capture the intuitive
notion of fairness.
\end{abstract}





\newpage

\tableofcontents

\newpage







\section{Introduction} 
\label{flowintro}

In optimization problems, 
a typical task is
to find an extreme element of a set $Q$ of `feasible' vectors, where
extreme means that we maximize (or minimize) a certain 
(linear or more general) objective function.  
A different (though related) concept in
optimization is when one is interested in finding an element of $Q$
whose components are distributed in a way which is felt 
the most uniform (fair, equitable, egalitarian).  
The term `fair' in the title of this paper refers to the intuitive meaning of the word.  
There may be various formal definitions for capturing this intuitive feeling.
For example, if the square-sum of the components is minimal, then the
distribution of the components is felt rather fair.  
Another possible way to formally capture fairness is to minimize the sum of the
absolute values of the pairwise differences of the components.  
A third possibility is lexicographic minimization.  
These definitions are equivalent in some cases while they are different in other situations.  
We should emphasize that the `fairness' concept shows
up in the literature in the most diverse contexts 
(such as fair resource allocation in operations research \cite{IK88,KSI13}, 
fair division of goods in economics \cite{Mou03,PR20}, 
load balancing in computer networks \cite{GGFS02,HLLT}, etc.).  
In the present work, however, fairness will be 
formulated into the concept of `decreasing minimality' (see, below).

An early example of a possible fairness concept is due to 
N. Megiddo \cite{Meg74,Meg77}, who introduced and solved
the problem of finding a (possibly fractional) maximum flow which is
\lq lexicographically optimal\rq \ on the set of edges 
leaving the source node.  
The problem, in equivalent terms, is as follows.  
Let $D=(V,A)$ be a digraph with a source-node $s$ and a sink-node $t$, and
let $S_{A}$ denote the set of edges leaving $s$.  
We assume that no edge enters $s$ and no edge leaves $t$.  
Let $g:A \rightarrow {\bf R}_{+}$ be
a non-negative capacity function on the edge-set.  
By the standard definition, an $st$-flow, or just a flow, is a function
$x:A \rightarrow {\bf R}_{+}$ 
for which $\varrho_{x}(v)=\delta_{x}(v)$
holds for every node $v\in V-\{s,t\}$. 
 (Here $\varrho_{x}(v):= \sum [x(uv):uv\in A]$ 
 and $\delta_{x}(v):= \sum [x(vu):  vu\in A]$.)  
The flow is called {\bf feasible} if $x\leq g$.  
The {\bf flow-amount} of $x$ is $\delta_{x}(s)$ 
which is equal to $\varrho_{x}(t)$.
We refer to a feasible flow with maximum flow-amount as a {\bf max-flow}.

Megiddo solved the problem of finding a feasible flow $x$ 
which is lexicographically optimal on $S_{A}$ in the sense 
that the smallest $x$-value on $S_{A}$ is as large as possible, 
within this, the second smallest (though not necessarily distinct) 
$x$-value on $S_{A}$ is as large as possible, and so on.  
It is a known fact
 (implied, for example, by the max-flow algorithm of Ford and Fulkerson
\cite{Ford-Fulkerson}) that a
lexicographically optimal flow is a max-flow.  
It is a basic property
of flows that for an integral capacity function $g$ there always
exists a max-flow which is integer-valued.  
On the other hand, an easy example \cite{FM21partA} shows that 
even when $g$ is integer-valued, the unique max-flow that
is lexicographically optimal on $S_{A}$ 
may not be integer-valued.

A member $x$ of a set $Q$ of vectors is called 
a {\bf decreasingly minimal} 
(dec-min,  for short) element of $Q$ 
if the largest component of $x$ is as small as possible, within this, 
the next largest (but not necessarily distinct)
component of $x$ is as small as possible, and so on.
The term `decreasing minimality' was introduced in \cite{FM21partA,FM21partB}
as one of the possible formulations of the intuitive notion of fairness.
Analogously, $x$ is an {\bf increasingly maximal} 
(inc-max) element of $Q$ 
if its smallest component is as large as possible, within this,
the next smallest component of $x$ is as large as possible, and so on.
Therefore increasing maximality is the same as Megiddo's lexicographic optimality
and `lexmin optimality' of Plaut and Roughgarden \cite{PR20},
whereas the notion of co-lexicographic optimality,  
introduced in
Fujishige \cite[page 264]{Fuj05book}, is the same as decreasing minimality.  
In general, a dec-min element is not necessarily inc-max,
and an inc-max element is not necessarily dec-min.
However, 
in Megiddo's problem where $Q$ is the restriction of a feasible 
maximum flow to $S_{A}$,
it is known that an element of $Q$ is dec-min if and only if it is inc-max.
Fujishige \cite{Fuj80,Fuj05book} proved that 
this equivalence is still true in a more general setting where
$Q$ is a base-polyhedron \cite{Fuj05book,Mdcasiam}.
He also proved that the (unique) dec-min element of $Q$ 
is the (unique) square-sum minimizer of $Q$.

In \cite{FM21partA} and \cite{FM21partB}, the present authors
solved the discrete counterpart of Megiddo's problem 
when the capacity function $g$ is integral and 
one is interested in finding an integral max-flow
whose restriction to the set $S_{A}$ of edges leaving $s$ is
increasingly maximal.  
This was actually a consequence of the more general result 
concerning dec-min elements of an M-convex set
(where an M-convex set \cite{Mdcasiam},
by definition, is the set of integral elements of an
integral base-polyhedron).
Among others, 
it was proved that an element $z$ is decreasingly minimal 
if and only if $z$ is increasingly maximal.  
It was also proved in \cite{FM21partA} that an element $z$ of an M-convex set 
is dec-min if and only if $z$ is square-sum minimizer.  
A strongly polynomial algorithm was also developed for finding a dec-min element.  
Since the restrictions of max-flows to
$S_{A}$ form a base-polyhedron, 
this gives an algorithm to
find an integral max-flow which is decreasingly minimal 
(and increasingly maximal) when restricted to $S_{A}$.

A closely related previous work 
is due to Kaibel, Onn, and Sarrabezolles \cite{KOS}.  
They considered (in an equivalent formulation) 
the problem of finding an integer-valued uncapacitated
$st$-flow with specified flow-amount $K$ 
which is decreasingly minimal on the whole edge-set $A$.  
They developed an algorithm which is polynomial in the size of 
digraph $D=(V,A)$ 
plus the value of $K$ but
is not polynomial
in the size of number $K$ 
(which is roughly $\lceil \log K\rceil $).
This is analogous to the well-known characteristic of the classic
Ford--Fulkerson max-flow algorithm \cite{Ford-Fulkerson},
where the running time is
proportional to the largest value $g_{\max}$ of the capacity
function $g$, and therefore this algorithm is not polynomial 
(unless $g_{\max}$ is small in the sense that 
it is bounded by a polynomial of $\vert A\vert $).  
It should also be mentioned that Kaibel et al.
considered exclusively 
 the uncapacitated $st$-flow problem, where no capacity
(upper-bound) restrictions are imposed on the edges.  
(For example, the flow-value on any edge is allowed to be $K$.)

In the present work, we consider the more general question when
$F\subseteq A$ is an arbitrarily specified subset of edges, and we are
interested in finding a feasible integral max-flow whose restriction
to $F$ is decreasingly minimal.  
This problem substantially differs from its special case 
with $F=S_{A}$ mentioned above  in that the set of
restrictions of max-flows to $F$ is not necessarily a base-polyhedron.
The significant difference is nicely demonstrated by the fact that an element $z$ 
of an M-convex set, as mentioned earlier, is dec-min if and only if 
it is inc-max if and only if it is a square-sum minimizer,
whereas 
these three criteria are (pairwise) different 
for integral feasible network flows
(see Section~\ref{SCconvmin}).
In this light, it is not surprising that the dec-min problem
for integral network flows is much
harder than for M-convex sets.

\begin{figure}\begin{center}
\includegraphics[height=35mm]{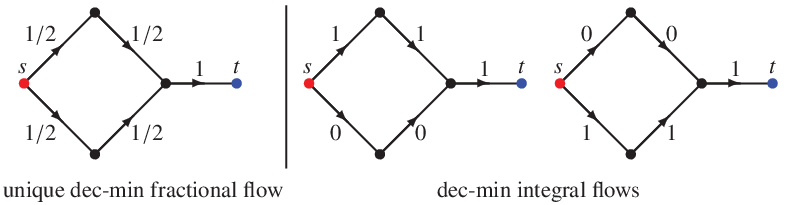}
\caption{Difference between fractional and integral dec-min flows}
\label{FGflowRZ}
\end{center}\end{figure}

We emphasize the fundamental difference between fractional and integral dec-min flows. 
Figure~\ref{FGflowRZ} demonstrates this difference for a simple example,
where all edges have a unit capacity ($g \equiv 1$)
and dec-min unit flows from $s$ to $t$ are considered for 
$F = A$ (all edges).
Whereas the dec-min fractional flow is uniquely determined,
there are two dec-min integral flows.

As the theory of network flows has a multitude of applications, the
algorithm presented in this paper may also be useful in these special cases.  
For example, the paper
by  Harvey, Ladner, Lov\'asz, and Tamir \cite{HLLT} 
considered the problem of finding a subgraph of a bipartite graph
$G=(S,T;E)$ for which
 the degree-sequence in $S$ is identically 1 and
the degree-sequence in $T$ is 
decreasingly minimal.  
This problem was extended to a more general setting 
(see \cite{FM21partB}) but the following
version needs the present general flow approach:
Find a subgraph of $G=(S,T;E)$ of $\gamma $ edges 
for which the degree-sequence on the whole node-set $S\cup T$ (or on an arbitrarily
specified subset of $S\cup T$) is decreasingly minimal.

Our main goal is to provide a description of the set of integral max-flows
which are dec-min on $F$ as well as a strongly polynomial algorithm to
find such a max-flow.  
The description makes it possible to solve
algorithmically even the minimum cost dec-min max-flow problem.
Instead of maximum $st$-flows, we consider the formally more general
(though equivalent) setting of modular flows which, however, allows a
technically simpler discussion.

It is quite natural to consider the dec-min problem over the intersection 
of two M-convex sets, which is called an M$_{2}$-convex set in the literature \cite{Mdcasiam}.
This problem is much harder than the dec-min problem over an M-convex set.  
The relationship of the difficulties is similar to 
that between the classic problems of finding a maximum weight basis of a matroid 
and finding a maximum weight common basis of two matroids 
(or more generally, between a maximum weight element of an M-convex set 
and of an M$_{2}$-convex set).  
An even more general framework is the set of integral submodular
flows, introduced by Edmonds and Giles \cite{EG77}, which
includes both standard integral network flows ($m$-flows) and M$_{2}$-convex sets.  
In \cite{FM22fairsbmflow}, we have worked out a strongly polynomial algorithm 
for finding an integral dec-min submodular flow.

The paper is organized as follows.  
In Section~\ref{SCdmfldef}, after introducing the basic definitions, 
we formulate Theorem~\ref{MAIN} which is the main theoretical result of the paper.  
This is proved in Section~\ref{SCdecmin} after the necessary structural results 
are developed in Sections \ref{flow2} and \ref{flow3}.
An important consequence of the characterization in Theorem~\ref{MAIN} 
is that it makes possible to manage algorithmically even 
the minimum edge-cost version of the integral dec-min flow problem.
Section~\ref{SCgall} provides an alternative characterization 
of $F$-dec-min integral feasible flows by developing extensions 
of such standard concepts 
from network optimization as improving di-circuits and feasible potentials.  
Section~\ref{veges} provides
a necessary and sufficient condition 
for the existence of an integral $F$-dec-min flow.  
Sections \ref{algo5}--\ref{SCalgsum}
are devoted to algorithmic aspects.
Sections \ref{algo5} and \ref{algo1}
describe strongly polynomial algorithms for each component, 
and Section~\ref{SCalgsum} shows how these components are synthesized.
Finally, in the supplementary 
Section~\ref{SCconcrem}
of the paper, we briefly outline two closely related topics:
fractional dec-min flows and the relation to convex minimization over flows.

\section{Decreasingly-minimal integer-valued feasible modular flows}
\label{SCdmfldef}

\subsection{Modular flows}

Let $D=(V,A)$ be a digraph endowed with integer-valued functions
$f:A \rightarrow {\bf Z} \cup \{-\infty \}$ 
and $g:A \rightarrow {\bf Z} \cup \{+\infty \}$  for which $f\leq g$.  
Here $f$ and $g$ are serving as lower and upper bound functions, respectively.  
An edge $e$ is called {\bf $(f,g)$-tight} or just {\bf tight} if $f(e)=g(e)$.  
The polyhedron $T(f,g):=\{x: f\leq x\leq g\}$ is called a {\bf box}.

We are given a finite integer-valued function $m$ on $V$ for which $\widetilde m(V)=0$.  
(Here and throughout,  $\widetilde m(X):=\sum [m(v):v\in X]$.)  
A {\bf modular flow} (with respect to $m$) or, for short, a {\bf mod-flow} $x$ 
is a finite-valued function on $A$ 
(or a vector in ${\bf R}\sp{A}$) 
for which $\varrho_{x}(v)-\delta_{x}(v)=m(v)$
for each node $v\in V$.  
When we want to emphasize the defining vector $m$, we speak of an {\bf $m$-flow}.

A mod-flow $x$ is called {\bf $(f,g)$-bounded} or {\bf feasible} if
$f\leq x\leq g$.  A circulation is an $m$-flow with respect to
$m\equiv 0$, and an $st$-flow of given flow-amount $K$ is also an
$m$-flow with respect to $m$ defined by
\begin{equation} 
m(v):= \begin{cases} 0 & \ \ \hbox{if}\ \ \ v\in V-\{s,t\}, 
\cr
 K & \ \ \hbox{if}\ \ \ v=t ,
\cr
 -K & \ \ \hbox{if}\ \ \ v=s. 
        \end{cases}
\end{equation}
Circulations form a subspace of ${\bf R}\sp{A}$ while the set of
mod-flows is an affine space.  
The set of feasible mod-flows, which is called a 
{\bf feasible mod-flow polyhedron,} may be viewed as the
intersection of this affine subspace with the box $T(f,g)$.  
It follows from this definition that the face of such a polyhedron is
also a feasible $m$-flow polyhedron.  
We note, however, that the projection along axes is not necessarily 
a feasible mod-flow polyhedron since its description may need 
an exponential number of inequalities 
while a feasible mod-flow polyhedron is described by at most 
$2\vert A\vert +\vert V\vert $ inequalities.

Let $Q=Q(f,g;m)$ denote the set of $(f,g)$-bounded $m$-flows.  
Hoffman's theorem \cite{Hoffman60} states that $Q$ is non-empty 
if and only if 
the Hoffman-condition 
$\varrho_{g}-\delta_{f}\geq \widetilde m$
holds, that is,
\begin{equation} 
\varrho_{g}(Z)-\delta_{f}(Z)\geq \widetilde m(Z) \quad
\hbox{for every}\ \ \ Z\subseteq V.
\label{(Hoffman)} 
\end{equation}
It is well-known that $Q$ is an integral polyhedron 
whenever $f$, $g$, and $m$ are integral vectors.  
In the integral case let $\odotZ{Q} = \odotZ{Q}(f,g;m)$ denote the set
of integral elements of $Q$, that is,
\begin{equation} \label{odottQ}
\odotZ{Q}:= Q \cap {\bf Z} \sp{A}. 
\end{equation}

In Section~\ref{flowintro} we introduced (the basic form of) the notion of decreasing minimality,
but we actually work with the following slightly extended definition.  
Let $F$ be a specified subset of $A$.  
We say that $z\in \odotZ{Q}(f,g;m)$ is 
{\bf decreasingly minimal on} $F$ (or {\bf $F$-dec-min} for short) 
if the restriction of $z$ to $F$ is decreasingly minimal
among the restrictions of the vectors in $\odotZ{Q}(f,g;m)$ to $F$.

Our first main goal is to prove the following characterization of the subset of
elements of $\odotZ{Q}$ which are decreasingly minimal on $F$.

\begin{theorem}  \label{MAIN} 
Let $D=(V,A)$ be a digraph endowed with integer-valued 
lower and upper bound functions 
$f:A \rightarrow {\bf Z} \cup \{-\infty \}$ 
and $g:A \rightarrow {\bf Z} \cup \{+\infty \}$ 
for which $f\leq g$.  
Let $m:V \rightarrow {\bf Z} $ be a function on $V$
with $\widetilde m(V)=0$ 
such that there exists an $(f,g)$-bounded $m$-flow.  
Let $F\subseteq A$ be a specified subset of edges such that
both $f$ and $g$ are finite-valued on $F$.  
There exists a pair $(f\sp{*},g\sp{*})$ of integer-valued functions on 
$A$ with $f\leq f\sp{*} \leq g\sp{*} \leq g$ 
(allowing $f\sp{*}(e)=-\infty $ and $g\sp{*}(e)=+\infty $ for $e\in A-F$) 
such that an integral $(f,g)$-bounded $m$-flow $z$ is
decreasingly minimal on $F$ if and only if $z$ is an integral 
$(f\sp{*},g\sp{*})$-bounded $m$-flow.  
Moreover, the box $T(f\sp{*},g\sp{*})$ is narrow 
on $F$ in the sense that $0\leq g\sp{*}(e)-f\sp{*}(e)\leq 1$ 
for every $e\in F$.  
\finbox
\end{theorem}

Our second main goal is to describe a strongly polynomial algorithm 
to compute $f\sp{*}$ and $g\sp{*}$.  
Once these bounds are available, one
is able to compute not only a single $(f,g)$-bounded integer-valued $m$-flow 
which is dec-min on $F$ but a minimum cost $F$-dec-min 
$m$-flow as well (with the help of a standard min-cost circulation algorithm).
Section~\ref{SCalgsum} 
summarizes what the various components of the whole algorithm aim at
and how these components are related to each other.

\begin{remark} \rm \label{RMnotused}
In Section~\ref{veges}, we shall consider the general
case when $f$ and $g$ are not required to be finite-valued on $F$.  
In this case, an $F$-dec-min $(f,g)$-feasible $m$-flow may not exist, 
and we shall provide a characterization for the existence.  
In Theorem~\ref{MAINb}, we shall show how Theorem~\ref{MAIN} 
can be extended to the case 
when only the existence of an $F$-dec-min $(f,g)$-feasible $m$-flow is assumed.
\finbox
\end{remark}

\begin{remark} \rm \label{RMincmaxflow}
One may also be interested in finding an (integral)
$(f,g)$-bounded $m$-flow $z$ which is {\bf increasingly maximal}
(inc-max) {\bf on} $F$ 
in the sense that the smallest $z$-value on $F$
is as large as possible, within this, the second smallest (but not
necessarily distinct) $z$-value on $F$ is as large as possible, and so on.  
(Megiddo \cite{Meg74}, \cite{Meg77}, for example,
considered the fractional inc-max problem for $st$-flows when $F$ was
the set of edges leaving $s$.)  
But an $(f,g)$-bounded $m$-flow $z$ is
increasingly maximal on $F$ precisely if $-z$ is a $(-g,-f)$-bounded
$(-m)$-flow which is dec-min on $F$, 
implying that the inc-max and the dec-min problems are equivalent for modular flows.  
Hence we concentrate throughout only on decreasing minimality.  
Note that in \cite{FM21partA}
we investigated these
problems for M-convex sets and proved that the two problems are not
only equivalent but they are one and the same in the sense that an
element $z$ of an M-convex set is dec-min if and only if $z$ is inc-max. 
 (As mentioned earlier, an M-convex set, by definition, is
nothing but the set of integral elements of an integral base-polyhedron).  
\finbox
\end{remark}

\begin{remark} \rm  \label{RMflowalg}
It is well-known that there are strongly polynomial algorithms that
find a feasible $m$-flow when it exists or find a subset $Z$ violating \eqref{(Hoffman)} 
(see, for example, appropriate variations of the
algorithms by Edmonds and Karp \cite{Edmonds-Karp}, Dinits \cite{Dinits}, 
or Goldberg and Tarjan \cite{Goldberg-Tarjan}).
Actually, when no feasible $m$-flow exists, not only a violating
subset can be computed but the most violating set as well, 
that is, a set $Z\sp{*}$ maximizing 
$\widetilde m(Z)-\varrho_{g}(Z)+\delta_{f}(Z)$.
Note that this latter function is fully supermodular
(see \eqref{supmoddef} for definition),
and there is a general algorithm to maximize an arbitrary supermodular function.  
The point here is that for finding $Z\sp{*}$ we do not have to rely on this
general algorithm since much simpler (and more efficient)
flow-techniques do the job.
\finbox
\end{remark}

\subsection{Approach of the proof of Theorem~\ref{MAIN}}
\label{approach}

By {\bf tightening an edge}
 $e$ we mean the operation that replaces the
bounding pair $(f(e),g(e))$ by 
\allowbreak
$(f'(e),g'(e))$ 
where 
$f(e)\leq f'(e)\leq g'(e)\leq g(e)$ and $g'(e)-f'(e)<g(e)-f(e)$.  
Note that tightening an edge does not necessarily make the edge tight.  
The approach of the proof is that we tighten edges as long as possible without
loosing any integral $m$-flow which is dec-min on $F$, and prove that
when no more tightening step is available for the current 
$(f\sp{*},g\sp{*})$ then 
every $(f\sp{*},g\sp{*})$-bounded integral $m$-flow is 
an $F$-dec-min element of $\odotZ{Q}(f,g;m)$.

A natural reduction step consists of removing a tight edge $e$ from $F$ 
(where $e$ could be tight originally or may have become tight
during a tightening step).  
This simply means that we replace $F$ by $F':=F-e$ 
(but keep $e$ in the digraph itself).  
Obviously, an $m$-flow $z$ is $F$-dec-min if and only if $z$ is $F'$-dec-min.
Therefore, we may always assume that $F$ contains no tight edges.

We say that an integral $(f,g)$-bounded $m$-flow $z$ is an $F$-{\bf max minimizer} 
if the largest component of $z$ in $F$ is as small as possible.  
Clearly, every $F$-dec-min $m$-flow $z\in \odotZ{Q}(f,g;m)$
is $F$-max minimizer.  
Let $\beta_{F}$ denote this smallest maximum value, that is,
\begin{equation} 
\beta_{F}:  = \min \{ \max \{z(a):  a\in F \}:  z\in \odotZ{Q}(f,g;m)\}.  
\label{(betaF)} 
\end{equation}
\noindent
Note that $\beta_{F}$ may be interpreted as the smallest integer 
for which there is an integer-valued feasible $m$-flow after
decreasing $g(e)$ to $\beta_{F}$ for each $e\in F$ with 
$g(e)>\beta_F$.  
In Section~\ref{algo1}, we shall describe how $\beta_{F}$ 
can be computed in strongly polynomial time with the help of
a discrete variant of 
the Newton--Dinkelbach algorithm and a standard max-flow algorithm, 
but for the proof of Theorem~\ref{MAIN} 
we assume that $\beta_{F}$ is available.  
Therefore, we can assume that $\max \{g(e):e\in F\} = \beta_{F}$ 
which is equivalent to requiring that $Q(f,g;m)$ is non-empty but 
$Q(f,g\sp{-};m)=\emptyset $ 
where $g\sp{-}$ arises from $g$ by subtracting 1 from 
$g(e)$ for each $e\in F$ with $g(e) = \beta_F$.

\section{Covering a supermodular function by a smallest subgraph}
\label{flow2}

A family of subsets is called {\bf laminar} if one of $X\subseteq Y$,
$Y\subseteq X$, $X\cap Y=\emptyset $ holds for every pair of its members.
We say that a digraph $D=(V,A)$ (or its edge-set $A$) {\bf covers} a
set-function $p$ if $\varrho_{D}(Z)\geq p(Z)$ for every subset $Z\subseteq V$,
where $\varrho_{D}$ is the in-degree function of $D$.
A set-function $p$ is called {\bf fully supermodular} or just 
{\bf supermodular} if
the supermodular inequality
\begin{equation} \label{supmoddef}
p(X)+p(Y) \leq p(X\cap Y) + p(X\cup Y) 
\end{equation}
holds for every pair of subsets $X$ and $Y$.  
When this inequality required only for intersecting pairs 
(that is, when $X\cap Y \ne \emptyset$), 
then we speak of an {\bf intersecting supermodular} function \cite{Frank-book}.

Let $p:2\sp{V} \rightarrow {\bf Z} \cup \{-\infty \}$ 
be an intersecting supermodular set-function on $V$ 
and let $D_{L}=(V,L)$ be a digraph covering $p$.  
We are interested in the minimum cardinality subset of edges of $D_{L}$ that covers $p$. 
Let $A_{L}$ denote the $(0,1)$-matrix
whose rows correspond to subsets $X$ of $V$ for which $p(X)>-\infty $
and the columns correspond to the edges in $L$.  
An entry of $A_{L}$ corresponding to $Z$ and $e$ is 1 
if $e$ enters $Z$ and 0 otherwise.
The following result was proved in \cite{FrankJ2} 
(see, also, Theorem 17.1.1 in the book \cite{Frank-book}).

\begin{theorem}  \label{kernel} 
Let $p$ be an intersecting supermodular set-function on $V$.  
The linear inequality system 
$[  A_{L}x_{L}\geq p, \  x_{L}\leq \underline{1}, \  x_{L}\geq 0  ]$ 
is totally dual integral (TDI). 
 (Hence) the primal linear program
\begin{equation} 
\min \{ \underline{1} x_{L}:  \ A_{L}x_{L}\geq p, \  x_{L}\leq \underline{1}, \  x_{L}\geq 0\}
\label{(primal)} 
\end{equation}
and the dual linear program 
\begin{equation} \max \{ yp - \underline{1} z:  yA_{L} - z \leq \underline{1}, \  (y,z)\geq 0 \} 
\label{(dual)} 
\end{equation}
have integer-valued optimal solutions,
where $\underline{1}$ denotes
the everywhere 1 vector of dimension $\vert L\vert $.
Moreover, there is an integer-valued dual optimum $(y\sp{*},z\sp{*})$ for which
its support family 
${\cal L}:=\{Z:  y\sp{*}(Z)>0\}$ is laminar.  
\finbox
\end{theorem}

For a family $\cal L$ of subsets of $V$,
let $\varrho_{L}({\cal L})$ denote
the number of edges entering at least one member of $\cal L$.  
The min-max theorem arising from Theorem~\ref{kernel} is as follows.

\begin{theorem}  \label{kernelminmax} 
Given a digraph $D_{L}=(V,L)$ covering an
intersecting supermodular function $p$, the minimum number of edges of
$D_{L}$ covering $p$ is equal to
\begin{equation} 
 \max \{ \varrho_{L}({\cal L}) \ 
 - \sum [\varrho_{L}(Z)-p(Z):  \ Z\in {\cal L}] \  \} ,
\label{(minmaxX)} 
\end{equation} 
where the maximum is taken over all laminar families $\cal L$ 
of subsets $Z$ of $V$ with $p(Z)>-\infty $. 
When $p$ is fully supermodular, the optimal laminar family ${\cal L}\sp{*}$ 
may be chosen as a chain of subsets 
$V_{1}\supset V_{2}\supset \cdots \supset V_{q}$ of $V$. 
\end{theorem}  

\Proof 
Suppose that we remove some edges from $L$ so that the set $X$
of the remaining edges continues to cover $p$.  
For each $Z\in {\cal L}$, 
the number of removed edges entering $Z$ is bounded by 
$\varrho_{L}(Z) - p(Z)$, 
and hence the number of removed edges entering at least one member 
of $\cal L$ is bounded from above by 
$\sum [\varrho_{L}(Z) - p(Z):  Z\in {\cal L}]$. 
On the other hand, the number of removed
edges entering at least one member of $\cal L$ is bounded from below by 
$\varrho_{L}({\cal L}) - \vert X\vert $. 
Therefore we have 
\[
\varrho_{L}({\cal L}) - \vert X\vert 
\leq \sum [\varrho_{L}(Z) - p(Z): Z\in {\cal L}], 
\] 
from which the trivial direction $\max \leq \min$ follows.

To see the reverse inequality, we have to find a covering 
$X\sp{*} \subseteq L$ of $p$ and a laminar family 
${\cal L}\sp{*}$ for which equality holds.  
To this end, let $x\sp{*}$ be a
$(0,1)$-valued optimal solution of the primal problem \eqref{(primal)}
in Theorem~\ref{kernel} and 
let $(y\sp{*},z\sp{*})$ be an integer-valued optimal solution 
of the dual problem for which its support family
${\cal L}\sp{*}$ is laminar.  
Then the subset $X\sp{*}:=\{e\in L :  x\sp{*}(e)=1\}$ 
is a smallest subset of $L$ covering $p$.

Observe that $y\sp{*}$ uniquely determines $z\sp{*}$, 
namely, $z\sp{*}(e)=0$ when $e$ enters no member of \ ${\cal L}\sp{*}$ and
\begin{equation} 
z\sp{*}(e) =
  \sum [ y\sp{*}(Z):  Z\in {\cal L}\sp{*}, e \ \hbox{enters}\ \ Z] \ -1 
\label{(zertek)} 
\end{equation}
when $e$ enters at least one member of ${\cal L}\sp{*}$.

\begin{claim} \label{CL01val}
The optimal $y\sp{*}$ may be chosen $(0,1)$-valued.
\end{claim}

\Proof 
Suppose that $(y\sp{*},z\sp{*})$ is an integer-valued dual optimum 
in which the sum of $y\sp{*}$-components is as small as possible.  
We show that $y\sp{*}$ is $(0,1)$-valued.  
Suppose indirectly that $y\sp{*}(Z)\geq 2$ for some set $Z$.  
In this case $z\sp{*}(e)\geq 1$ for every edge $e$ entering $Z$.  
If we decrease $y\sp{*}(Z)$ by 1 and decrease $z\sp{*}(e)$ by 1 
on every edge $e$ entering $Z$, then the resulting $(y',z')$ 
is also a dual feasible solution for which
\[
y\sp{*}p - \underline{1} z\sp{*} 
 \ \geq \  y'p - \underline{1} z' 
 \ = \  y\sp{*}p - \underline{1} z\sp{*} - p(Z) + \varrho_{L}(Z) 
 \  \geq \  y\sp{*}p - \underline{1} z\sp{*} ,
\] 
where the last inequality follows from the assumption that $D_{L}$ covers $p$ 
and hence $\varrho_{L}(Z)\geq p(Z)$.  
Therefore we have equality throughout and 
hence $(y',z')$ is also an optimal dual solution, contradicting the
minimal choice of $y\sp{*}$.  
Thus Claim \ref{CL01val} is proved.
\finbox 
\medskip

By the claim, \eqref{(zertek)} simplifies as follows:
\begin{equation} z\sp{*}(e) 
= \hbox{$[$the number of members of $\cal L$ entered by $e] \ -1.$}\ 
\label{(zertekb)} 
\end{equation}
\noindent
Now the dual optimum value is:
\begin{eqnarray} 
&& y\sp{*}p - \underline{1} z\sp{*}
\nonumber\\
 &&= \sum [p(Z): Z\in {\cal L}\sp{*}] 
 -  \sum [z\sp{*}(e) :  e\in L \ \hbox{enters a
member of}\ {\cal L}\sp{*} ] 
\nonumber\\
 &&= \sum [p(Z):  Z\in {\cal L}\sp{*}] 
\nonumber\\ 
&& \ \ {}- \sum [( \hbox{the number of members of}\ {\cal L}\sp{*} 
\hbox{entered by}\ e) - 1:  \ e \ \hbox{enters a member of}\ {\cal L}\sp{*}] 
\nonumber\\ 
&&= \sum [p(Z):  Z\in {\cal L}\sp{*}] 
- \sum [\varrho_{L}(Z) :  Z\in {\cal L}\sp{*}]
  + \varrho_{L}({\cal L}\sp{*}) 
\nonumber\\ 
&&= \varrho_{L}({\cal L}\sp{*}) \ 
- \sum [\varrho_{L}(Z)-p(Z):  \ Z\in {\cal L}\sp{*}] .  
\label{(dualopt)}
\end{eqnarray}
\noindent
Therefore $\vert X\sp{*}\vert $ is equal to the value in
\eqref{(dualopt)}, from which the non-trivial direction $\max \geq \min $ follows, 
implying the requested $\min = \max$.

To see the last statement of the theorem, consider 
an optimal laminar family $\cal L$ with a minimum number of members.  
We claim that $\cal L$ is a chain of subsets when $p$ is fully supermodular.  
Suppose, indirectly, that $\cal L$ has two disjoint members and let $X$ and $Y$
be disjoint members of $\cal L$ whose union is maximal.  
Then the family ${\cal L}'$ obtained from $\cal L$ by replacing $X$ and $Y$
with their union $X\cup Y$ is also laminar.  
By the full supermodularity of $p$, 
we have \ 
$\sum [p(Z):  \ Z\in {\cal L}] \leq \sum [p(Z):  \ Z\in {\cal L}']$.  
Furthermore, 
\[
\varrho_{L}({\cal L})
- \sum [\varrho_{L}(Z):  \ Z\in {\cal L}] = \varrho_{L}({\cal L}') -
\sum [\varrho_{L}(Z):  \ Z\in {\cal L'}].  
\]
\noindent
Therefore ${\cal L}'$ is also a dual optimal laminar family,
contradicting the minimal choice of $\cal L$. 
This completes the proof of Theorem~\ref{kernelminmax}.
\finbox \finboxHere

\begin{theorem}   \label{optkritX} 
Let $D_{L}=(V,L)$ be a digraph covering a
fully supermodular function $p$.  There is a chain ${\cal C}\sp{*}$ of
subsets $V_{1}\supset V_{2}\supset \cdots \supset V_{q}$ of $V$ with
$p(V_{i})> -\infty $ such that a subset $X\subseteq L$ is a minimum
cardinality subset of edges covering $p$ if and only if the following
three optimality criteria hold.

\noindent 
{\rm (A)} \ For every $V_{i}$, \ $\varrho_{X}(V_{i})= p(V_{i})$.  \

\noindent
{\rm (B)} \ Every edge in $X$ enters at least one $V_{i}$.
(Equivalently, if $e\in L$ enters no $V_{i}$, then $e\not \in X$.)

\noindent
{\rm (C)} \ Every edge in $L-X$ enters at most one $V_{i}$.
(Equivalently, if $e\in L$ enters at least two $V_{i}$'s, then $e\in X$.)  
\end{theorem}

\Proof 
Let ${\cal C}\sp{*}$ denote the optimal chain of subsets
$V_{1}\supset V_{2}\supset \cdots \supset V_{q}$ \ 
given in Theorem~\ref{kernelminmax}.  
This corresponded to a special integer-valued
solution $(y\sp{*},z\sp{*})$ to the dual linear program \eqref{(dual)}
where $y\sp{*}$ was actually $(0,1)$-valued and $y\sp{*}$ 
(or its support family ${\cal C}\sp{*}$) determined uniquely $z\sp{*}$.  
Namely, $z\sp{*}(e)$ was 0 when $e$ did not enter any $V_{i}$, 
and $z\sp{*}(e)$ was the number of $V_{i}$'s entered by $e$ minus 1 
when $e$ entered at least one $V_{i}$.

Since both the primal and the dual variables in the linear programs in
Theorem~\ref{kernel} are non-negative, the optimality criteria ($=$
complementary slackness conditions) of linear programming require that
if a primal variable is positive, then the corresponding dual
inequality holds with equality, and symmetrically, if a dual variable
is positive, then the corresponding primal inequality holds with equality.

Let $x\sp{*}$ be a $(0,1)$-valued primal solution and let 
$X\sp{*}:=\{e\in L:  x\sp{*}(e)=1\}$ 
be the corresponding set of edges that covers $p$.  
The optimality criterion concerning the dual variable $y\sp{*}$, 
requires that if $y\sp{*}(Z)=1$ 
(that is, if $Z$ is one of the sets $V_{i}$), 
then the corresponding primal inequality holds with equality.  
That is, 
$\varrho_{X\sp{*}}(V_{i})= \varrho_{x\sp{*}}(V_{i}) = p(V_{i})$, 
which is just Criterion (A).

The optimality criterion concerning the primal variable $x\sp{*}$
requires that if $x\sp{*}(e)=1$ for an edge $e$ 
(that is, if $e\in X\sp{*}$), 
then the corresponding dual inequality holds with equality.
Hence $e$ must enter at least one $V_{i}$ (as $z\sp{*}(e)\geq 0$), which
is just Criterion (B).

Finally, the optimality criterion concerning 
the dual variable $z\sp{*}(e)$ requires that if $z\sp{*}(e)>0$ 
(that is, if $e$ enters at least two $V_{i}$'s), 
then the corresponding primal inequality is met by equality, 
that is, $x\sp{*}(e)=1$ or equivalently $e\in X\sp{*}$, which
is just Criterion (C).
\finbox

\section{$L$-upper-minimal $m$-flows} 
\label{flow3}

Let $D=(V,A)$ be a digraph and $m:V \rightarrow {\bf Z}$ a function
with $\widetilde m(V)=0$.  
Let $f:A \rightarrow {\bf Z} \cup \{-\infty \}$ 
and $g:A \rightarrow {\bf Z} \cup \{+\infty \}$ 
be bounding functions with $f\leq g$.  
Let $L$ be a subset of $A$ for which $-\infty <f(e)<g(e)<+\infty $ 
for every $e\in L$.  
(That is, $f(e)$ may be $-\infty $ and 
$g(e)$ may be $+\infty $ only if $e\in A-L$.)
We say that an $(f,g)$-bounded integer-valued $m$-flow $x$ is 
{\bf $L$-upper-minimal} or that $x$ is an {\bf $L$-upper-minimizer} 
if the number of $g$-saturated edges in $L$ is as small as possible, 
where an edge $e\in L$ is called {\bf $g$-saturated} if $x(e)=g(e)$.  
In this section, we are interested in characterizing the $L$-upper-minimizer
integral $(f,g)$-bounded $m$-flows.  
For the proof of Theorem~\ref{MAIN}, however, 
we will use this characterization only in the special case when 
$L:=\{e:  e\in F, g(e)=\beta_{F}\}$, 
that is, $g(e)$ is the same value for each element $e$ of $L$.  
The only reason for this more general setting is to get a clearer 
picture of the background.

Let $g\sp{-}:=g-\chi_{L}$, that is,
\begin{equation} 
g\sp{-}(e):= \begin{cases} 
  g(e)-1 & \ \ \hbox{if \ \ $e\in L$, }
  \cr 
  g(e) & \ \ \hbox{if \ \ $e \in A-L$.  } 
   \end{cases} 
\end{equation}
\noindent
Since $g(e)<+\infty $ for $e\in L$, $g\sp{-}\not =g$.  
By the hypothesis, $L$ contains no tight edges 
and hence $f\leq g\sp{-}$. 
Define a set-function $p$ as follows:
\begin{equation} 
p:= \widetilde m - \varrho_{g\sp{-}} + \delta_{f}.  
\label{(pdef)}
\end{equation}
Since $g\sp{-}\geq f$, the function $\varrho_{g\sp{-}}-\delta_{f}$ 
is fully submodular and hence $p$ is fully supermodular.  
Furthermore,
$p(Z)>-\infty $ precisely if $\varrho_{g}(Z) - \delta_{f}(Z)<+\infty $.

The following lemma states a basic fact,
which will be used several times in the proofs of
Theorems \ref{minL} and \ref{optkritx}.

\begin{lemma} \label{LMxX} 
{\rm (A)} 
If $x$ is an integer-valued $(f,g)$-bounded
$m$-flow, and $X\subseteq L$ is the set of $g$-saturated $L$-edges,
(that is, $X:=\{e\in L:  x(e)=g(e)\})$, then $X$ covers $p$.  
\ 
{\rm (B)} 
If a subset $X\subseteq L$ covers $p$, then there is an integer-valued
$m$-flow which is $(f,g\sp{-}+\chi_X)$-bounded.
\end{lemma}  

\Proof 
(A) For every subset $Z\subseteq V$, we have
\[
  \widetilde m(Z) =  \varrho_{x}(Z) - \delta_{x}(Z) 
  \leq [\varrho_{g\sp{-}}(Z) + \varrho_{X}(Z)] - \delta_{f}(Z),
\]
from which
\[
 \varrho_{X}(Z) \geq \widetilde m(Z) - \varrho_{g\sp{-}}(Z) + \delta_f(Z) =p(Z) ,
\]
as required.

(B) \ It follows from the hypothesis 
$\varrho_{X}\geq p = \widetilde m - \varrho_{g\sp{-}} + \delta_{f}$ 
that 
$\varrho_{g\sp{-}} + \varrho_{X} - \delta_{f}\geq \widetilde m$.  
Then Hoffman's theorem implies that
there is an integer-valued 
$(f,g\sp{-}+\chi_X)$-bounded $m$-flow.
\finbox

\medskip

The next lemma shows a key fact 
which connects
an $L$-upper-minimizer flow to the general framework of supermodular covering
presented in Section~\ref{flow2}.

\begin{lemma}  \label{LMminmin} 
An integer-valued $(f,g)$-bounded $m$-flow $x$ is an $L$-upper-minimizer 
if and only if 
$X:=\{e\in L:  x(e)=g(e)\}$ \
is a smallest subset of $L$ covering $p$.  
\end{lemma}

\Proof
The proof consists of the following two claims.
\begin{claim} \label{Xfromx} 
If $x$ is an $L$-upper-minimizer $(f,g)$-bounded
$m$-flow, then $X:=\{e\in L:  x(e)=g(e)\}$ is a smallest subset of $L$
covering $p$.
\end{claim}

\Proof 
By Part (A) of Lemma~\ref{LMxX}, we know that $X$ covers $p$.
Let $X'\subseteq L$ be an arbitrary cover of $p$, that is,
\[
\varrho_{X'} \geq \widetilde m -\varrho_{g\sp{-}} +\delta_{f},
\]
or equivalently,
\[
 \varrho_{X'}+ \varrho_{g\sp{-}} -\delta_{f} \geq \widetilde m. 
\]

By Part (B) of Lemma~\ref{LMxX},
there exists an integer-valued $m$-flow
$x'$ which is $(f,g\sp{-}+\chi_{X'})$-bounded.  
Hence every $g$-saturated $L$-edge (with respect to $x'$) belongs to $X'$.  
Since $x$ is an $L$-upper-minimizer, 
it follows that $\vert X\vert \leq \vert X'\vert $, 
that is, $X$ is indeed a smallest subset of $L$ covering $p$.
\finbox

\begin{claim} \label{smallcover} 
If $X\sp{*} \subseteq L$ is a smallest subset of $L$ covering $p$, 
then every integer-valued $(f,g\sp{-}+\chi_{X\sp{*}})$-bounded $m$-flow $x\sp{*}$ 
is an $L$-upper-minimizer $(f,g)$-bounded $m$-flow.
\end{claim}

\Proof 
Let $X':=\{e\in L:  x\sp{*}(e)=g(e)\}$.  
By Lemma~\ref{LMxX},
$X'$ covers $p$ and hence 
$\vert X\sp{*}\vert \leq \vert X'\vert $. 
Since $x\sp{*}$ is $(f,g\sp{-}+\chi_{X\sp{*}})$-bounded, 
it follows that $x\sp{*}$ admits at most 
$\vert X\sp{*}\vert $ $g$-saturated $L$-edges
from which 
$\vert X\sp{*}\vert \geq \vert X'\vert $. 
Therefore $\vert X\sp{*}\vert =\vert X'\vert $ 
and thus $x\sp{*}$ saturates a minimum number of elements of $L$, 
that is, $x\sp{*}$ is an $L$-upper-minimizer.
\finbox 

\medskip

This completes the proof of Lemma~\ref{LMminmin}.  
\finbox \finboxHere

\medskip

The following min-max theorem shall be the basis of an optimality criterion 
for the decreasing minimality of a feasible  $m$-flow,   
and hence it serves as a stopping rule of the algorithm.
In the following two theorems, we use notations 
$D=(V,A), L, f,g, m$ introduced in the first paragraph of this section.  
In particular, we assume that $-\infty < f(e)<g(e)<+\infty$ holds for every edge $e\in L$.

\begin{theorem}  \label{minL} 
The minimum number of $g$-saturated $L$-edges in
an $(f,g)$-bounded integer-valued $m$-flow is equal to
\begin{equation} 
 \max \{ \varrho_{L}({\cal C})
  \ - \sum [\varrho_{g}(Z)- \delta_{f}(Z) - \widetilde m(Z):  \ Z\in {\cal C}]  \} ,
 \label{(minmaxL)} 
\end{equation} 
where the maximum is taken over all chains $\cal C$ of subsets $Z$ of $V$
with $\varrho_{g}(Z)- \delta_{f}(Z)< +\infty $, 
and $\varrho_{L}({\cal C})$ 
denotes the number of $L$-edges entering at least one member of ${\cal C}$.  
In particular, if the minimum is zero, the maximum is
attained at the empty chain.  
\end{theorem}

\Proof 
Let $x$ be an $(f,g)$-bounded integer-valued $m$-flow with a minimum number of
$g$-saturated $L$-edges.  
Let $X=\{e\in L:  x(e)=g(e)\}$, 
that is, $X$ is the set of $g$-saturated $L$-edges.  
By Lemma~\ref{LMminmin}, 
$X$ is a smallest subset of $L$ covering $p$.

Apply Theorem~\ref{kernelminmax} to the digraph $D_{L}=(V,L)$ and to the
set-function $p$ defined in \eqref{(pdef)}.  
In this case, $p$ is fully supermodular from which we obtain that 
\begin{align*}
\vert X\vert &= \max \{
 \varrho_{L}({\cal C}) \ - \sum [\varrho_{L}(Z)-p(Z):  \ Z\in {\cal C}] :
 \ {\cal C} \ \hbox{a chain of subsets of}\ \ V\} 
\\ & =
\max \{ \varrho_{L}({\cal C}) \ - 
  \sum [\varrho_{g}(Z)-\delta_{f}(Z) - \widetilde m(Z) : \ Z\in {\cal C}] : 
   \ {\cal C} \ \hbox{a chain of subsets of} \ V \}, 
\end{align*}
as required.
This completes the proof of Theorem~\ref{minL}.
\finbox

\medskip 
Our next goal is to obtain optimality criteria for $L$-upper-minimizer $m$-flows.

\begin{theorem}   \label{optkritx} 
There is a chain ${\cal C}\sp{*}$ of subsets 
$V_{1}\supset V_{2}\supset \cdots \supset V_{q}$
of $V$ with $\varrho_{g}(V_{i})-\delta_{f}(V_{i})<+\infty $ such that an
integer-valued $(f,g)$-bounded $m$-flow $z$ is an $L$-upper-minimizer
if and only if the following optimality criteria hold.

\noindent 
{\rm (O1)} \ $z(e)=f(e)$ for every edge $e\in A$ leaving a set $V_{i}$, \

\noindent 
{\rm (O2)} \ $z(e)=g(e)$ for every edge $e\in A-L$ entering a set $V_{i}$,

\noindent 
{\rm (O3)} \ $g(e)-1\leq z(e)\leq g(e)$ for every edge $e\in L$
entering exactly one $V_{i}$,

\noindent 
{\rm (O4)} \ $z(e)=g(e)$ for every edge $e\in L$ entering at least two $V_{i}$'s,

\noindent 
{\rm (O5)} \ $f(e)\leq z(e)\leq g(e)-1$ for every edge $e\in L$
neither entering nor leaving any $V_{i}$.
\end{theorem}

\Proof
We shall apply Theorem~\ref{optkritX} to the digraph $D_{L}=(V,L)$ and to
the set-function $p$ defined in \eqref{(pdef)}, 
and consider the chain 
${\cal C}\sp{*}=\{V_{1},\dots ,V_{q}\}$ 
ensured by the theorem, where
$V_{1}\supset \cdots \supset V_{q}$.  
Since $p(V_{i})$ is finite for each $i=1,\dots ,q$, 
so is $\varrho_{g}(V_{i})-\delta_{f}(V_{i})$.  
Note that both $f(e)$ and $g(e)$ are finite for each edge $e\in L$ 
and for each edge leaving or entering a member of ${\cal C}\sp{*}$.

To see the necessity of the conditions (O1)--(O5),
suppose that $x\sp{*}$ is an
integer-valued $(f,g)$-bounded $m$-flow which is an $L$-upper-minimizer.  
By Lemma~\ref{LMminmin},
the set $X\sp{*}:=\{e\in L:  x\sp{*}(e)=g(e)\}$
is a smallest subset of $L$ covering $p$.  
Hence the optimality criteria (A), (B), and (C) in Theorem~\ref{optkritX} hold.

By Property (A), 
$\varrho_{X\sp{*}}(V_{i})= p(V_{i})$ 
for every $V_{i}$, which is equivalent to
\begin{equation} 
\varrho_{g\sp{-}}(V_{i}) + \varrho_{X\sp{*}}(V_{i}) -\delta_{f}(V_{i}) 
= \widetilde m(V_{i}) , 
\end{equation}
from which
\[
 \widetilde m(V_{i}) = \varrho_{x\sp{*}}(V_{i}) - \delta_{x\sp{*}}(V_{i})
 \leq \varrho_{g\sp{-}}(V_{i}) + \varrho_{X\sp{*}}(V_{i}) - \delta_{f}(V_{i})
= \widetilde m(V_{i}).
\]
Hence we have equality throughout, in particular,
\begin{equation} \varrho_{x\sp{*}}(V_{i})
 = \varrho_{g\sp{-}}(V_{i}) + \varrho_{X\sp{*}}(V_{i}) 
 \quad [ = \widetilde m(V_{i}) + \delta_{f}(V_{i})] 
\label{(xopt1)}
\end{equation} and 
\begin{equation} 
  \delta_{x\sp{*}}(V_{i}) = \delta_{f}(V_{i}).  
 \label{(xopt2)}
\end{equation}

The equality in \eqref{(xopt2)} shows that (O1) holds.  
Condition \eqref{(xopt1)} implies for an edge $e\in A-L$ entering a $V_{i}$ 
that $x\sp{*}(e)=g\sp{-}(e) = g(e)$ and hence (O2) holds.  
Condition \eqref{(xopt1)} implies for an edge $e\in L$ entering a $V_{i}$ that
$g(e)-1\leq x\sp{*}(e)\leq g(e)$ 
and hence (O3) holds.
By Property (C), if an edge $e\in L$ enters at least two $V_{i}$'s, then
$e\in X\sp{*}$ and hence $x\sp{*}(e)=g(e)$, that is, (O4) holds.
To see (O5), let $e\in L$ be an edge neither entering nor leaving any $V_{i}$.  
By Property (B), $e\not \in X\sp{*}$ and 
hence $x\sp{*}(e)\leq g(e)-1$, from which (O5) follows.

To see the sufficiency of the conditions (O1)--(O5),
let $z$ be an integer-valued
$(f,g)$-bounded $m$-flow satisfying the five conditions in the theorem.  
Let $X:=\{e\in L:  z(e)=g(e)\}$.  
By Part (A) of Lemma~\ref{LMxX}, $X$ covers $p$.  
We claim that $X$ meets the three optimality criteria in Theorem~\ref{optkritX}.  
Let $V_{i}$ be a member of chain ${\cal C}\sp{*}$.

(O2) implies that 
\[ 
\sum [z(e):  e \in A-L, e \ \hbox{enters}\ V_{i}] =
\sum [g(e):  e \in A-L, e \ \hbox{enters}\ V_{i}].
\]
\noindent
 From the definition of $X$, we have 
\[ 
\sum [z(e):  e \in X, e \ \hbox{enters}\ V_{i}] 
= \sum [g(e):  e \in X, e \ \hbox{enters}\ V_{i}].
\]
(O3) implies that 
\[ 
\sum [z(e):  e \in L-X, e \
\hbox{enters}\ V_{i}] = \sum [g(e)-1:  e \in L-X, e \ \hbox{enters}\ V_{i}].  
\]
\noindent
By merging these three equalities, we obtain
\[
\varrho_{z}(V_{i}) = \varrho_{g\sp{-}}(V_{i}) + \varrho_{X}(V_{i}).
\]
\noindent
Furthermore, (O1) implies that 
\[
\delta_{z}(V_{i})=\delta_{f}(V_{i}) ,
\] 
from which
\[
\widetilde m(V_{i}) = \varrho_{z}(V_{i}) - \delta_{z}(V_{i}) 
= \varrho_{g\sp{-}}(V_{i}) + \varrho_{X}(V_{i}) -\delta_{f}(V_{i}), 
\] 
that is,
\[
\varrho_{X}(V_{i}) 
= \widetilde m(V_{i}) - \varrho_{g\sp{-}}(V_{i}) +\delta_{f}(V_{i})=p(V_{i}),
\]
showing that Property (A) in Theorem~\ref{optkritX} holds indeed.

To see Property (B), let $e\in X$ ($\subseteq L$) be an edge.  
Then $z(e)=g(e)$ and, by (O5), $e$ enters or leaves a $V_{i}$.  
But $e$ cannot leave any $V_{i}$ since if it did, then (O1) 
would imply $z(e)=f(e)$ and this would contradict the assumption 
that $L$ contains no tight edge.  
Therefore $e$ must enter a $V_{i}$, that is, (B) holds indeed.

To see Property (C), let $e$ be an edge in $L$ which enters at least
two $V_{i}$'s.  By (O4), $z(e)=g(e)$ and hence $e\in X$, that is, (C) holds.

By Theorem~\ref{optkritX}, $X$ is a smallest subset of $L$ covering $p$.  
By Lemma~\ref{LMminmin}, $x$ is an $L$-upper-minimizer
$(f,g)$-bounded $m$-flow, as stated in the theorem.
This completes the proof of Theorem~\ref{optkritx}.
\finbox

\medskip

In Section~\ref{algo5}, we describe an algorithmic proof of Theorem~\ref{optkritx}.  
The algorithm will compute in strongly polynomial
time an $(f,g)$-bounded $L$-upper-minimizer integral $m$-flow along
with the optimal chain described in the theorem.

\section{Description of $F$-dec-min $m$-flows:  Proof of Theorem~\ref{MAIN}} 
\label{SCdecmin}

After preparations in Sections \ref{flow2} and \ref{flow3}, 
we turn to our main goal of proving Theorem~\ref{MAIN}.  
As before, let $D=(V,A)$ be a digraph and $F\subseteq A$ a specified subset of edges.  
We assume that the underlying undirected graph of $D$ is connected.  
Let $f:A \rightarrow {\bf Z} \cup \{-\infty \}$ 
and $g:A \rightarrow {\bf Z} \cup \{+\infty \}$ 
be bounding functions with $f\leq g$.  
We require $-\infty <f(e)\leq g(e)<+\infty $ for every $e\in F$.  
Let $m:V \rightarrow {\bf Z} $ 
be a function on the node-set for which
there is an integer-valued $(f,g)$-bounded $m$-flow 
(that is, $\widetilde m(V)=0$ and Hoffman's condition \eqref{(Hoffman)} holds).
Recall from \eqref{odottQ}
that $\odotZ{Q}= \odotZ{Q}(f,g;m)$ denotes the set of
integer-valued $(f,g)$-bounded $m$-flows.

In the proof we shall use induction on $\vert F\vert $. 
Since $f\sp{*}:=f$ and $g\sp{*}:=g$ 
clearly meet the requirements of the theorem when
$F=\emptyset $, we can assume that $F$ is non-empty.  We observed
already in Section~\ref{approach} that it suffices to prove Theorem~\ref{MAIN} 
in the special case when $F$ contains no tight edge,
therefore we assume throughout that $f(e)<g(e)$ for each edge $e\in F$.

Let $\beta =\beta_{F}$ 
denote the smallest integer for which $\odotZ{Q}$
has an element $z$ satisfying $z(e)\leq \beta$ for every edge $e\in F$
(cf., \eqref{(betaF)}).
In Section~\ref{algo1}, we shall work out an algorithm to compute
$\beta_{F}$ in strongly polynomial time.  Since we are interested in
$F$-dec-min members of $\odotZ{Q}$, we may assume that the largest
$g$-value of the edges in $F$ is this $\beta$. 
Let $L:=\{e\in F: g(e)=\beta \}$.  
Now Hoffman's condition \eqref{(Hoffman)} holds but,
since $F$ contains no tight edges and since $\beta$ is minimal, after
decreasing the $g$-value of the elements of $L$ from $\beta$ to $\beta -1$, 
the resulting function $g\sp{-}:=g-\chi_{L}$ violates \eqref{(Hoffman)}, 
that is, $Q(f,g\sp{-};m)=\emptyset $. 
Summing up, we shall rely on the following notation and assumptions:
\begin{equation} 
\begin{cases} & 
\hbox{$F$ is non-empty and contains no $(f,g)$-tight edges,}
\cr & 
\hbox{$\beta :=\max \{g(e):  e\in F\}$},  
\cr & 
\hbox{$L:=\{e\in F:  g(e) =\beta \}$},
\cr & 
\hbox{$g\sp{-}:=g-\chi_{L}$}, 
\cr & 
\hbox{$\odotZ{Q} = \odotZ{Q}(f,g;m)$ is non-empty,
}\ 
\cr & 
\hbox{$\odotZ{Q}(f,g\sp{-};m)$ is empty}.  
\end{cases}
\label{(hypo)} 
\end{equation}

\medskip

As a preparation for deriving the main result Theorem~\ref{MAIN}, 
we need the following relaxation of decreasing minimality.  
We call a member $z$ of $\odotZ{Q}$ {\bf pre-decreasingly minimal} 
({\bf pre-dec-min}, for short) {\bf on} $F$ 
if the number $\mu $ of edges $e$ in $L$ with 
$z(e) = \beta$ is as small as possible.  
Obviously,
if $z$ is $F$-dec-min, then $z$ is pre-dec-min on $F$.  
By applying Theorem~\ref{optkritx} to the present special case, 
we obtain the following characterization of pre-dec-min elements.

\begin{theorem}   \label{Fdecmin} 
Given \eqref{(hypo)}, there is a chain
${\cal C}'$ of non-empty proper subsets 
$V_{1}\supset V_{2}\supset \cdots \supset V_{q}$ 
of $V$ with 
$\varrho_{g}(V_{i})-\delta_{f}(V_{i})<+\infty $
such that a member $z$ of $\odotZ{Q}$ is pre-dec-min on $F$ 
if and only if 
the following optimality criteria hold:  
\medskip

\noindent 
{\rm (O1)} \ $z(e)=f(e)$ for every edge $e\in A$ leaving a member
of ${\cal C}'$, \ \smallskip

\noindent 
{\rm (O2)} \ $z(e)=g(e)$ for every edge $e\in A-L$ entering a
member of ${\cal C}'$, \smallskip

\noindent 
{\rm (O3)} \ $\beta -1\leq z(e)\leq \beta$ for every edge $e\in
L$ entering exactly one member of ${\cal C}'$, \smallskip

\noindent
{\rm  (O4)} \ $z(e)=\beta$ for every edge $e\in L$ entering at
least two members of ${\cal C}'$, \smallskip

\noindent 
{\rm  (O5)}  \ $f(e)\leq z(e)\leq \beta-1$ for every edge $e\in L$
neither entering nor leaving any member of ${\cal C}'$.  
\finbox 
\end{theorem}

We call a chain with the properties in the theorem a {\bf dual optimal chain}.  
Section~\ref{algo5} describes a strongly polynomial algorithm
for computing a dual optimal chain.

Define the bounding pair $(f'(e),g'(e))$ for each edge $e$, as follows.  
For $e\in L$, let
\begin{equation} 
(f'(e),g'(e)):= 
\begin{cases}
 (\beta ,\beta ) & \hbox{if $e$ enters at least two members of ${\cal C}'$}, 
  \cr 
(\beta -1,\beta ) & \hbox{if $e$ enters exactly one member of ${\cal C}'$}, 
  \cr
(f(e),f(e)) & \hbox{if $e$ leaves a member of ${\cal C}'$}, 
 \cr
(f(e),\beta -1) & \hbox{if $e$ neither leaves nor enters any member of ${\cal C}'$}. 
\end{cases} 
\label{(f'g'def1)} 
\end{equation}
\noindent
For $e\in A-L$, let 
\begin{equation} (f'(e),g'(e)):= 
\begin{cases}
(g(e),g(e)) & \hbox{if $e$ enters a member of ${\cal C}'$}, 
\cr
(f(e),f(e)) & \hbox{if $e$ leaves a member of ${\cal C}'$}, 
\cr
(f(e),g(e)) & \hbox{if $e$ neither leaves nor enters any member of ${\cal C}'$}. 
\end{cases} 
\label{(f'g'def2)} 
\end{equation}
\noindent
It follows from this definition that $f\leq f'\leq g'\leq g$. 
 Let
\begin{equation} 
 \hbox{ $\odotZ{Q'}:  = \odotZ{Q}(f',g';m)$.  }\ 
\label{(Q')} 
\end{equation}

\begin{lemma}  \label{f'g'} 
{\rm (A)}
 \ An $m$-flow $z\in \odotZ{Q}$ is pre-dec-min on $F$ 
if and only if $z\in \odotZ{Q'}$.  
\ {\rm (B)} \ An $m$-flow $z\in \odotZ{Q}$ is $F$-dec-min 
if and only if $z$ is an $F$-dec-min element of $\odotZ{Q'}$.  
\end{lemma}

\Proof 
Theorem~\ref{Fdecmin} immediately implies the equivalence in Part (A).  
To see Part (B), suppose first that $z$ is an $F$-dec-min element of $\odotZ{Q}$.  
Then $z$ is surely $F$-pre-dec-min in $\odotZ{Q}$
and hence, by Part (A), $z$ is in $\odotZ{Q'}$.  
If, indirectly, $\odotZ{Q'}$ had an element $z'$ 
which is decreasingly smaller on $F$ than $z$, 
then $z$ could not have been an $F$-dec-min element of $\odotZ{Q}$.  
Conversely, let $z'$ be an $F$-dec-min element of 
$\odotZ{Q'}$ and suppose indirectly that $z'$ 
is not an $F$-dec-min element of $\odotZ{Q}$.  
Then any $F$-dec-min element $z$ of $\odotZ{Q}$ 
is decreasingly smaller on $F$ than $z'$.  
But any $F$-dec-min element of $\odotZ{Q}$ is pre-dec-min on $F$ and hence, 
by Part (A), $z$ is in $\odotZ{Q'}$, contradicting the assumption that $z'$ was an
$F$-dec-min element of $\odotZ{Q'}$.  
\finbox 
\medskip

Theorem~\ref{MAIN} will be an immediate consequence of the following result.

\begin{theorem}  \label{reduction} 
Given \eqref{(hypo)}, 
there is a pair $(f',g')$ of integer-valued functions on $A$ with 
$f\leq f'\leq g'\leq g$ and a set $F'\subset F$ 
such that an element $z$ of $\odotZ{Q}$
is an $F$-dec-min member of $\odotZ{Q}$ 
if and only if $z$ is an
$F'$-dec-min member of
 \ $\odotZ{Q'}=\odotZ{Q}(f',g';m)$.  
In addition,
the box $T(f',g')$ is narrow on $F-F'$ in the sense that 
$0\leq g'(e)-f'(e)\leq 1$ holds for every $e\in F-F'$.  
\end{theorem}

\Proof 
Let ${\cal C}'$ be the chain ensured by Theorem~\ref{Fdecmin},
let $(f',g')$ be the pair of bounding functions defined in
\eqref{(f'g'def1)} and \eqref{(f'g'def2)}, and let 
$\odotZ{Q'}:  = \odotZ{Q}(f',g';m)$.
Let $L'$ denote the subset of $L$ consisting of those
elements of $L$ that enter at least one member of ${\cal C}'$.

\begin{claim} \label{claimL'}
The set $L'\subseteq L$ is non-empty.  
\end{claim}

\Proof 
Let $z$ be an element of $\odotZ{Q}$ which is pre-dec-min on $F$.  
By Part (A) of Lemma~\ref{f'g'}, $z\in \odotZ{Q'}$.  
By \eqref{(hypo)},
there is an edge $e$ in $F$ for which $z(e)=\beta =g(e)$, and hence $e\in L$.  
Since $g(e)=z(e)\leq g'(e)\leq g(e)$ and $F$ contains no
$(f,g)$-tight edges, we have $f(e) < g(e) =g'(e) =\beta$. 
This and definition \eqref{(f'g'def1)} imply that $e$ enters at least one member
of ${\cal C}'$.  
\finbox \medskip

Since $L'\not =\emptyset $ by the claim, we have
\[ 
\hbox{ $F':=F-L'$ is a proper subset of $F$}. 
\]
\noindent
We are going to show that $(f',g')$ and $F'$ meet the requirements of
the theorem.  
Call two vectors in ${\bf Z}\sp{A}$ 
{\bf value-equivalent}
on $L'$ if their restrictions to $L'$ (that is, their projection 
to ${\bf Z}\sp{L'}$), when both arranged in a decreasing order, are equal.

\begin{lemma}  \label{valeq} 
The members of $\odotZ{Q'}$ are value-equivalent on $L'$.  
\end{lemma}

\Proof 
By Part (A) of Lemma~\ref{f'g'}, the members of $\odotZ{Q'}$
are exactly those elements of $\odotZ{Q}$ which are pre-dec-min on $F$.
Hence each member $z$ of $\odotZ{Q'}$ has the same number $\mu $ of
edges $e$ in $L$ with $z(e)=\beta$.

As $F$ contains no $(f,g)$-tight edges, we have $z(e)\leq g'(e)\leq \beta -1$ 
for every edge $e\in L-L'$ and hence each element $e$ of $L$
with $z(e)=\beta$ belongs to $L'$, from which
\[
 \vert \{e\in L':  z(e)= \beta \}\vert = \mu .
\]
\noindent
Furthermore, we have $f'(e)\geq \beta -1$ for every element $e$ of
$L'$, from which $L'$ has exactly $\vert L'\vert -\mu $ edges with
$z(e)=\beta -1$, implying that the members of $\odotZ{Q'}$ are indeed
value-equivalent on $L'$.  
\finbox \medskip

Part (B) of Lemma~\ref{f'g'} implies that the $F$-dec-min elements of
$\odotZ{Q}$ are exactly the $F$-dec-min elements of $\odotZ{Q'}$, 
and hence it suffices to prove that an element $z$ of $\odotZ{Q'}$ 
is an $F$-dec-min member of $\odotZ{Q'}$
if and only if 
$z$ is an $F'$-dec-min member of \ $\odotZ{Q'}$.  
But this latter equivalence is
an immediate consequence of Lemma~\ref{valeq}.

To prove the last part of Theorem~\ref{reduction}, recall that
$F-F'=L'$ and $L'$ consisted of those elements of $L$ that enter at
least one member of ${\cal C}'$.  
But the definition of $(f',g')$ in \eqref{(f'g'def1)} 
implies that $\beta -1\leq f'(e)\leq g'(e)=\beta$
for every element $e$ of $L'$, that is, the box $T(f',g')$ is indeed
narrow on $F-F'$.  
This completes the proof of Theorem~\ref{reduction}.
\finbox \finboxHere 
\medskip

\paragraph{Proof of Theorem~\ref{MAIN}}  
\ We use induction on $\vert F\vert $. 
Since $f\sp{*}:=f$ and $g\sp{*}:=g$ clearly meet the
requirements of the theorem when $F=\emptyset $, we can assume that
$F$ is non-empty.  
As before, we may assume that $F$ contains no
$(f,g)$-tight edges.  
By Theorem~\ref{reduction}, it suffices to prove
Theorem~\ref{MAIN} for $\odotZ{Q}(f',g';m)$ and $F'$.  
But this follows
by induction since $F'$ is a proper subset of $F$.  
\finbox

\paragraph{Cheapest integral $F$-dec-min $m$-flows} \ 
In Sections \ref{algo5} and \ref{algo1}, 
we shall describe 
a strongly polynomial algorithm to compute
$(f\sp{*},g\sp{*})$ in Theorem~\ref{MAIN}.  
Once these bounding
functions are available, we can immediately solve the problem of
computing a cheapest integral $F$-dec-min $(f,g)$-bounded $m$-flow
with respect to a cost-function $c:A \rightarrow {\bf R}$.
By Theorem~\ref{MAIN}, 
this latter problem is nothing but a minimum cost 
$(f\sp{*},g\sp{*})$-bounded $m$-flow problem, which can indeed be solved by a
minimum cost feasible circulation algorithm.  
In the literature there
are several strongly polynomial algorithms for the cheapest
circulation problem, the first one was due to Tardos \cite{Tardos1}.

\section{Characterization by improving di-circuits and by feasible potential-vectors} 
\label{SCgall}

Let $D=(V,A)$, $F$, $f$, $g$, $m$ be the same as in Theorem~\ref{MAIN}.  
Let $\odotZ{Q} = \odotZ{Q}(f,g;m)$ denote the set 
of integral $(f,g)$-bounded $m$-flows.  
We assume that $\odotZ{Q}$ is non-empty but the properties
in \eqref{(hypo)} are not a priori expected.  
For an element $z\in \odotZ{Q}$, 
let $D_{z}=(V,A_{z})$ denote the standard auxiliary digraph
associated with $z$, that is,
\[
 A_{z}:=\{uv:  uv\in A, z(uv)<g(uv)\} \cup \{vu :  uv\in A, z(uv)>f(uv)\}.
\]
\noindent
An edge $uv\in A_{z}$ is called a forward edge when
$z(uv)<g(uv)$ and a backward edge when $z(vu)>f(vu)$.

Theorem~\ref{MAIN} provided a characterization for the set of
$F$-dec-min elements of $\odotZ{Q}$, namely, an element $z\in \odotZ{Q}$
is $F$-dec-min precisely if $f\sp{*} \leq z \leq g\sp{*}$.  
The goal of
this section is to describe a different characterization for 
$z\in \odotZ{Q}$ to be decreasingly minimal on $F$, consisting of two
equivalent properties.  
(For a comparison of the previous and this new
characterizations, see Remark \ref{char.link}.)  
For the first one, we introduce a simple and natural way 
to obtain from $z$ a decreasingly smaller feasible $m$-flow 
by improving $z$ along an appropriate di-circuit of $D_{z}$.  
For the second property, by extending the
standard notion of feasible potentials, we introduce feasible potential-vectors.  
The main result of the section
(Theorem~\ref{gall.main} in Section~\ref{SCcharthm})  states (roughly)
that the following three properties for $z$ are pairwise equivalent:
\ (A) \ $z$ is dec-min on $F$, 
\ (B) \ no di-circuit improving $z$ exists, 
and 
\ (C) \ there exists a feasible potential-vector.

\subsection{Feasible potential-vectors}

Let $c:A_{0} \rightarrow {\bf R} $ be a cost-function defined on the
edge-set of a digraph $D_{0}=(V,A_{0})$.  
A di-circuit $C$ of $D_{0}$ is
called negative (with respect to $c$) if the total $c$-cost
$\widetilde c(C) =\sum [c(e):e\in C]$ of $C$ is negative.  
In the literature, $c$ is called {\bf conservative} 
if $D_{0}$ admits no negative di-circuit.  
A function $\pi :V \rightarrow {\bf R} $ is called 
a {\bf $c$-feasible potential} 
if $\pi (v)-\pi (u)\leq c(uv)$ holds for every edge $uv$ of $D_{0}$.  
A classic result of Gallai is as follows.

\begin{theorem}[Gallai]  \label{Gallai} 
Given a digraph $D_{0}=(V,A_{0})$ and a
cost-function $c:  A_{0}  \rightarrow {\bf R}$, there exists a
$c$-feasible potential $\pi :  V  \rightarrow {\bf R}$ 
if and only if
$c$ is conservative.  If $c$ is conservative and integer-valued, 
then $\pi $ can be chosen integer-valued, as well.  
\finbox
\end{theorem}

Given two $k$-dimensional vectors $\underline{x}=(x_{1},x_{2},\dots ,x_{k})$ and
$\underline{y}=(y_{1},y_{2},\dots ,y_{k})$, we say that $\underline{x}$ 
is {\bf lexicographically smaller} than $\underline{y}$, 
in notation $\underline{x}\prec \underline{y}$, 
if $\underline{x}\not = \underline{y}$ and $ x_{i}< y_{i}$
 where $i$ denotes the first component in which they differ.  
We write $\underline{x}\preceq \underline{y}$ 
if $\underline{x}= \underline{y}$ or $\underline{x}\prec \underline{y}$.  
Note that the relation $\preceq $ is a total ordering of the elements of ${\bf R}\sp{k}$.

Let $\underline{c} :  A_{0}  \rightarrow {\bf R}\sp{k}$ 
be a vector-valued function on the edge-set of $D_{0}=(V,A_{0})$ 
that assigns a vector $\underline{c}(e)=(c_{1}(e), c_{2}(e),\dots ,c_{k}(e))$ 
to each edge $e$ of $D_{0}$.  
We call a vector-valued function $\underline{\pi} :  V \rightarrow {\bf R}\sp{k}$
on the node-set $V$ {\bf $\underline{c}$-feasible} or just feasible if 
\begin{equation} 
\underline{\pi} (v) - \underline{\pi} (u) \preceq \underline{c}(uv)
\label{feasiblepotvec}
\end{equation} 
holds for every edge $uv$ of $D_{0}$.

A di-circuit $C$ is said to be {\bf $\underline{c}$-negative} if the sum
$\widetilde {\underline{c}}(C)
=(\widetilde c_{1}(C),\widetilde c_{2}(C),\dots ,\widetilde c_{k}(C))$ 
of the $\underline{c}$-vectors assigned to its edges is
lexicographically smaller than the $k$-dimensional zero vector $\underline{0}_{k}$.  
The vector-valued function $\underline{c}$ is {\bf conservative} if
$D_{0}$ has no $\underline{c}$-negative di-circuit.

The following Gallai-type theorem specializes to Theorem~\ref{Gallai}
in case $k=1$, but in its proof we rely on Theorem~\ref{Gallai}.

\begin{theorem}  \label{ext-Gallai} 
Given a digraph $D_{0}=(V,A_{0})$ and a
vector-valued function $\underline{c}:  A_{0} \rightarrow {\bf R}\sp{k}$ on its
edge-set, there exists a $\underline{c}$-feasible potential-vector 
$\underline{\pi} : V \rightarrow {\bf R}\sp{k}$ 
if and only if 
$\underline{c}$ is conservative,
that is, $D_{0}$ admits no $\underline{c}$-negative di-circuit.  
If $\underline{c}$ is integer vector-valued and conservative, then a $\underline{c}$-feasible 
$\underline{\pi} $ can be chosen to be integer vector-valued.
\end{theorem}

\Proof 
Let $C$ be a di-circuit of $D_{0}$ whose nodes, in cyclic order,
are $v_{1},v_{2},\dots ,v_{q}$.  
Accordingly, the edges of $C$ are
$e_{1}=v_{1}v_{2}, e_{2}=v_{2}v_{3}, \dots , e_{q}=v_{q}v_{1}$.  
Let $\underline{\pi} $ be a $\underline{c}$-feasible potential-vector. 
Then
\begin{align*}
 \underline{0}_{k} & = [\underline{\pi} (v_{2}) - \underline{\pi} (v_{1})] 
 + [\underline{\pi} (v_{3}) - \underline{\pi} (v_{2})] + \cdots 
 + [\underline{\pi} (v_{1}) - \underline{\pi} (v_{q})] 
\\ &
\preceq  \sum [\underline{c}(e_{i}):  i=1,\dots ,q] = \widetilde {\underline{c}}(C).
\end{align*}

To see the reverse direction, we apply induction on $k$.  
When $k=1$, we are back at Theorem~\ref{Gallai}.  
Suppose now that $k\geq 2$, and
assume that $D_{0}$ admits no $\underline{c}$-negative di-circuit.

Consider the functions $c_{i}:A_{0} \rightarrow {\bf R} $ 
formed by the $i$-th components of $\underline{c}$ \ ($i=1,\dots ,k$).  
As $\underline{c}$ is conservative, so is $c_{1}$, 
that is $\widetilde c_{1}(C)\geq 0$ for every di-circuit $C$.  
By Theorem~\ref{Gallai}, there exists a $c_{1}$-feasible potential 
$\pi_{1}:  V \rightarrow {\bf R} $ (which is
integer-valued when $c_{1}$ is integer-valued).  
Let $A_{1}$ denote 
the following set of edges:
\[
  A_{1}:= \{uv\in A_{0}:  \pi_{1}(v)-\pi_{1}(u)=c_{1}(uv)\}.
\]
\noindent
Let $k':=k-1$ and $\underline{c}':=(c_{2},c_{3},\dots ,c_{k})$.  Then
$\underline{c}'$ is conservative in $D_{1}=(V,A_{1})$ since $\underline{c}$ is
conservative and $\pi_{1}(v)-\pi_{1}(u)=c_{1}(uv)$ holds for every edge
$uv$ in $A_{1}$.  
By induction, there is a $(k-1)$-dimensional potential-vector, 
$\underline{\pi }'= (\pi_{2},\dots ,\pi_{k})$
\ which is $\underline{c}'$-feasible on the edges in $A_{1}$.  
Let $\underline{\pi }:  = (\pi_{1},\pi_{2},\dots ,\pi_{k})$.  
Then $\underline{\pi }$ is $\underline{c}$-feasible on the edges in $A_{1}$.  
Moreover, $\pi_{1}(v)-\pi_{1}(u)<c_{1}(uv)$ for every edge $uv\in A_{0}-A_{1}$, 
and hence $\underline{\pi }$ is $\underline{c}$-feasible on
these edges, as well. 
 \finbox \medskip

\begin{remark} \rm \label{RMpotveccomp1}
A standard result of network flow theory is that
if the cost-function $c$ 
in Theorem~\ref{Gallai} is conservative, then a
$c$-feasible potential $\pi $ can be computed in polynomial time with
the help of the Bellman--Ford algorithm
(see, e.g., \cite[page 108]{Sch03}).  
Because the proof of Theorem~\ref{ext-Gallai} 
applies Theorem~\ref{Gallai} iteratively $k$ times, 
we can conclude that 
if the cost-vector $\underline{c}$ in the theorem is conservative 
and $k$ is polynomially bounded by $|A_{0}|$, 
then a  $\underline{c}$-feasible potential-vector $\underline{\pi}$ 
can be computed in polynomial time, and this $\underline{\pi}$ 
is an integral vector when $\underline{c}$ is an integral vector.  
We note that Theorem~\ref{ext-Gallai} 
and this algorithmic approach 
will be applied in the proof of 
Theorem~\ref{gall.main} where $k\leq 2 |F|$. 
\finbox
\end{remark}

\subsection{Improving di-circuits}

Let $A_{+}$ and $A_{-}$ be two disjoint sets and let $A_{*}:=A_{+}\cup A_{-}$.
Let $x$ be an integer-valued function on $A_{*}$.  
As a preparatory lemma, we develop an equivalent condition  for the function 
\begin{equation} 
x':= x+ \chi_{A_{+}}- \chi_{A_{-}} 
\label{(x'def)} 
\end{equation} 
to be decreasingly smaller than $x$.  
To this end, define $x\sp{*}:  A_{*} \rightarrow {\bf Z} $, as follows:
\begin{equation} 
x\sp{*}:= x - \chi_{A_{-}}.  
\label{(x*def1)} 
\end{equation}
Let $\lambda_{1}>\lambda_{2}>\cdots >\lambda_{h}$ 
denote the distinct values of the components of $x\sp{*}$.
We assign a $h$-dimensional vector $\underline{c}'(e)$ to every element $e\in A_{*}$, 
as follows:
\begin{equation}
 \underline{c}'(e):= \begin{cases}
  \   \underline{\varepsilon}'_{i} & \hbox{if $e\in A_{+}$ and $x\sp{*}(e)=\lambda_{i}$, }\ 
\cr 
  -\underline{\varepsilon}'_{i} &  \hbox{if $e\in A_{-}$ and $x\sp{*}(e)=\lambda_{i}$, }\ 
   \end{cases}
\label{(c'def)} 
\end{equation}
where $\underline{\varepsilon}'_{i}$ is the $h$-dimensional unit
vector $(0,\dots ,0,1,0,\dots ,0)$ whose $i$-th component is 1.

\begin{lemma} \label{LMekvi.0} 
$x' <_{\rm dec} x$ \ if and only if \
$\widetilde {{\underline{c}}'} (A_{*}) \prec \underline{0}_{h}$.  
\end{lemma}  

\Proof 
Induction on $\vert A_{*}\vert$. 
If $\vert A_{*}\vert =0$, then the statement of the lemma is void, 
so suppose that $A_{*}\not =\emptyset $. 
If $A_{-}=\emptyset$ and $A_{+}\not =\emptyset$, 
then $x' >_{\rm dec} x$ and $\widetilde {{\underline{c}}'} (A_{*}) \succ \underline{0}_{h}$, 
and hence neither of the two inequalities in the lemma holds.  
If $A_{-}\not =\emptyset$ and $A_{+}=\emptyset $, 
then $x' <_{\rm dec} x$ and $\widetilde {{\underline{c}}'} (A_{*}) \prec \underline{0}_{h}$, 
and hence both of the two inequalities in the lemma hold.  
So we can suppose that $A_{-}\not =\emptyset $ and $A_{+}\not =\emptyset $.

Let $e_{+}$ be an element of $A_{+}$ 
for which $\lambda_{i}=x\sp{*}(e_{+})$ is maximum, 
and let $e_{-}$ be an element of $A_{-}$ for which 
$\lambda_{j}=x\sp{*}(e_{-})$ 
is maximum.  
If $\lambda_{i} > \lambda_{j}$, then $x' >_{\rm dec} x$ 
and $\widetilde {{\underline{c}}'} (A_{*}) \succ \underline{0}_{h}$,
and hence neither of the two inequalities in the lemma holds.  
If $\lambda_{i} < \lambda_{j}$, 
then $x' <_{\rm dec} x$ and $\widetilde {{\underline{c}}'} (A_{*}) \prec \underline{0}_{h}$,
that is, both of the inequalities in the lemma hold.

In the remaining case, when $\lambda_{i}=\lambda_{j}$, we have $x(e_{+})+1=x(e_{-})$.  
Define $A'_{+}:=A_{+}-e_{+}$, $A'_{-}:=A_{-} - e_{-}$, and let
$A'_{*}:=A_{*}-\{e_{-},e_{+}\}$.  
Observe that the restriction of $x'$ to $A'_{*}$ 
is decreasingly smaller than the restriction of $x$ to $A'_{*}$
precisely if $x' <_{\rm dec} x$.  On the other hand, 
$\widetilde {{\underline{c}}'} (A'_{*}) = \widetilde {{\underline{c}}'} (A_{*})$ 
and hence 
$\widetilde {{\underline{c}}'} (A'_{*}) \prec \underline{0}_{h}$ 
precisely if $\widetilde {{\underline{c}}'} (A_{*}) \prec \underline{0}_{h}$.  
Since $\vert A'_{*}\vert < \vert A_{*}\vert $, 
we are done by induction.  
\finbox

\medskip

After this preparation, we return to $D=(V,A)$
with $F\subseteq A$ and $z\in \odotZ{Q}=\odotZ{Q}(f,g;m)$.  
Let $D_{z}=(V,A_{z})$ be the auxiliary digraph associated with $z$. 
We call a di-circuit $C$ of $D_{z}$ {\bf $z$-improving on $F$}
 (or just $z$-improving) if 
$z'\in \odotZ{Q}$ is decreasingly smaller than 
$z$ on $F$, where $z'(uv)$ is defined for $uv\in A$, as follows:
\begin{equation} 
z'(uv):= \begin{cases} 
 z(uv)+1 & \hbox{if $uv$ is a forward edge of $C$}, 
\cr
 z(uv)-1 & \hbox{if $vu$ is a backward edge of $C$}, 
\cr
z(uv) & \hbox{otherwise}. 
\end{cases} 
\label{(f'g'def22)} 
\end{equation}
\noindent
Note that the definition of $D_{z}$ implies that $z'$ is
indeed in $\odotZ{Q}$.

Let $F_{z}$ denote the subset of $A_{z}$ corresponding to $F$ 
(that is,
for $uv\in F$, if $z(uv)<g(uv)$, then the forward edge $uv$ belongs to $F_{z}$, 
while if $z(uv)>f(uv)$, then the backward edge $vu$ belongs to $F_{z}$).  
The sets of forward and backward edges in $F_{z}$ are denoted
by $F_{\bf f}$ and $F_{\bf b}$, respectively.  
(The subscripts ${\bf f}$ and ${\bf b}$ refer to {\bf f}orward and {\bf b}ackward.)

Define a function $z\sp{*}$ on $F_{z}$, as follows:
\begin{equation} 
z\sp{*}(uv):= \begin{cases} 
z(uv) & \hbox{if $uv\in F_{\bf f}$},
\cr 
z(vu)-1 & \hbox{if $uv\in F_{\bf b}$}.  
\end{cases}
\label{(z*def)} 
\end{equation}
\noindent
Let $\gamma_{1} > \gamma_{2}>\cdots >\gamma_{k}$ 
denote the distinct values of $z\sp{*}$,
where $k \leq 2|F|$.
Let $\underline{\varepsilon}_{i}$ denote the $k$-dimensional unit-vector
$(0,\dots ,0,1,0, \dots ,0)$ 
whose $i$-th component is 1. We assign a $k$-dimensional vector
$\underline{c}(e)$ to every edge $e$ of $D_{z}$, as follows:
\begin{equation} 
\underline{c}(e):= \begin{cases} 
\ \underline{0}_{k} & \hbox{if $e\in A_{z}-F_{z}$}, 
\cr
\ \underline{\varepsilon }_{i} & \hbox{if $e\in F_{\bf f}$ and $z\sp{*}(e)=\gamma_{i}$, }\ 
\cr 
\ -\underline{\varepsilon }_{i} & \hbox{if $e\in F_{\bf b}$ and $z\sp{*}(e)=\gamma_{i}$}. 
\end{cases} 
\label{(cdef)} 
\end{equation}

\begin{lemma}  \label{LMekvi}
 A di-circuit $C$ of $D_{z}$ is $z$-improving on $F$
if and only if \ $\widetilde {\underline{c}}(C)\prec \underline{0}_{k}$.  
\end{lemma}

\Proof 
Let $A_{+}:=\{uv:  uv \in F_{\bf f}\cap C\}$, 
$A_{-}:=\{uv:  vu\in F_{\bf b}\cap C\}$, 
and $A_{*}:=A_{+}\cup A_{-}$.  
Note that $A_{*} \subseteq A$.  
Let $x$ denote the restriction of $z$ to $A_{*}$.  
Then $x'$ defined in \eqref{(x'def)} is the restriction of $z'$ to $A_{*}$, 
and $x\sp{*}$ defined in \eqref{(x*def1)} 
is the restriction of $z\sp{*}$ to $A_{*}$.  
Let $\lambda_{1}>\lambda_{2}>\cdots >\lambda_{h}$ 
denote the distinct values of $x\sp{*}$, 
and consider the vector $\underline{c}'$ defined in \eqref{(c'def)}.  
Note that $\{\lambda_{1},\lambda_{2},\dots ,\lambda_{h}\}$ 
is a subsequence of $\{\gamma_{1},\gamma_{2},\dots ,\gamma_{k}\}$,
in particular, $h\leq k$.  
Observe that $C$ is $z$-improving if and only if
 $x'$ is decreasingly smaller than $x$.  
Also observe that
$\widetilde {\underline{c}}(C) \prec \underline{0}_{k}$ if and only if 
$\widetilde {{\underline{c}}'}  (A_{*}) \prec \underline{0}_{h}$.  
Then we are done by Lemma~\ref{LMekvi.0}.
\finbox

\subsection{Minimizing the number of saturated edges}
\label{SCminsatedge}

Let $\beta :=\max \{g(e):  e\in F\}$ and 
let $L:=\{e\in F:  g(e)=\beta \}$.  
We assume that $-\infty < f(e) < \beta$ for every edge $e\in L$,
while $f(e)=-\infty$ and $g(e)=+\infty$ are allowed for edges $e$ in $A-L$.  
The goal of this section
is to characterize $(f,g)$-bounded integral $m$-flows 
which saturate a minimum number of $L$-edges.

We need the following standard characterization of cheapest feasible $m$-flows.

\begin{lemma} \label{javkor} 
Let $D_{1}=(V,A_{1})$ be a digraph endowed 
with a cost function $c_{1}:  A_{1}\rightarrow {\bf R} $ 
and a pair $(f_{1},g_{1})$ of bounding-functions on $A_{1}$.  
For an $(f_{1},g_{1})$-bounded integral $m$-flow $x$, 
let $D_{x}=(V,A_{x})$ denote the auxiliary digraph, endowed
with a cost-function $c_{x}:A_{x}\rightarrow {\bf R}$, 
in which $uv\in A_{x}$ is a forward edge 
if $x(uv)<g_{1}(uv)$, for which $c_{x}(uv):=c_{1}(uv)$, 
and $vu\in A_{x}$ is a backward edge 
if $x(uv)>f_{1}(uv)$, for which $c_{x}(vu):=-c_{1}(uv)$.  
Then $x$ is a cheapest $(f_{1},g_{1})$-bounded integral $m$-flow 
if and only if there is no negative di-circuit in $D_{x}$ 
(or in other words, $c_{x}$ is conservative).  
\finbox 
\end{lemma}

In order to characterize integral $(f,g)$-bounded $m$-flows for which
the number of $g$-saturated (that is, $\beta$-valued) edges in $L$ is minimum, 
we introduce a parallel copy $e'$ of each $e\in L$.  Let $L'$
denote the set of new edges.  
Let $A_{1}:=A\cup L'$ and $D_{1}:=(V,A_{1})$.
Define $g\sp{-}$ on $A$ by $g\sp{-}:=g-\chi_{L}$, that is, we reduce
$g(e)$ from $\beta$ to $\beta -1$ for each $e\in L$.

Let $f_{1}$ and $g_{1}$ be bounding functions on $A_{1}$ defined by
\begin{equation}   \label{f1g1def}
 f_{1}(e) := \begin{cases}
           f(e) & \hbox{if $e\in A$}, \cr 
           0 & \hbox{if $e\in L'$}, 
            \end{cases} 
\qquad
 g_{1}(e) := \begin{cases} 
      g\sp{-}(e) & \hbox{if $e\in A$}, \cr 
          1 & \hbox{if $e\in L'$}.
            \end{cases} 
\end{equation}
\noindent
Let $c_{1}$ be a $(0,1)$-valued cost-function on $A_{1}$ defined by
\begin{equation}   \label{c1def}
 c_{1}(e):= \begin{cases}
 0 & \hbox{if $e\in A$}, \cr 
  1 & \hbox{if $e\in L'$}.
        \end{cases} 
\end{equation}

\begin{lemma} \label{minbeta}
\quad
\\ 
{\rm (A)} \ 
If $z$ is an integral $(f,g)$-bounded $m$-flow in $D$ 
having $\mu $ edges in $L$ with $z(e)=\beta$, 
then there exists an integral $(f_{1},g_{1})$-bounded $m$-flow $z_{1}$ in $D_{1}$
for which $c_{1}z_{1}=\mu$. 
\\
{\rm (B)} \ 
If $z_{1}$ is a minimum $c_{1}$-cost integer-valued $(f_{1},g_{1})$-bounded $m$-flow in $D_{1}$, 
then there is an
$(f,g)$-bounded $m$-flow $z$ in $D$ for which the number of 
edges in $L$ with $z(e)=\beta$ is $c_{1}z_{1}$.  
\end{lemma}

\Proof (A) Let $z$ be an $m$-flow given in Part (A), and let
$X:=\{e\in L:  z(e)=\beta \}$.  
Let $X'$ denote the subset of $L'$ corresponding to $X$.  
Define an $m$-flow $z_{1}$ in $D_{1}$ as follows:
\begin{equation}  \label{(z1def)} 
z_{1}(e):= \begin{cases}
      z(e) & \hbox{if $e\in A-X$}, \cr 
            \beta -1 &  \hbox{if $e\in X$}, \cr 
           1 & \hbox{if $e\in X'$}, \cr 
           0 & \hbox{if $e\in L'-X'$}. 
\end{cases} 
\end{equation}
\noindent
Then $z_{1}$ is an $(f_{1},g_{1})$-bounded $m$-flow in $D_{1}$ whose
$c_{1}$-cost is $\vert X\vert =\mu$.

(B) Let $z_{1}$ be an $m$-flow given in Part (B) of the lemma.  
Observe that if $z_{1}(e')=1$ for some $e'\in L'$, 
then $z_{1}(e)=g_{1}(e)=\beta -1$
where $e$ is the edge in $L$ corresponding to $e'$.  
Indeed, if we had $z_{1}(e)\leq \beta -2$, 
then the $m$-flow obtained from $z_{1}$ by adding $1$ to $z_{1}(e)$ 
and subtracting 1 from $z_{1}(e')$ would be of smaller cost.  
It follows that the $m$-flow $z$ in $D$ defined by
\begin{equation} 
 z(e):= \begin{cases} 
 z_{1}(e)+z_{1}(e') & \hbox{if $e\in L$}, \cr
  z_{1}(e) & \hbox{if $e\in A-L$ }
        \end{cases} 
\label{(zdef)} 
\end{equation}
\noindent 
is an $(f,g)$-bounded $m$-flow in $D$, for which the number
of $\beta$-valued $L$-edges is exactly the $c_{1}$-cost of $z_{1}$.  
\finbox \medskip

\begin{corollary} \label{ekviz-z1} 
An integral $(f,g)$-bounded $m$-flow $z$ in $D$
with $\max \{z(e):e\in L\}\leq \beta$ minimizes 
the number of the $\beta$-valued edges in $L$ 
if and only if 
the $(f_{1},g_{1})$-bounded $m$-flow $z_{1}$ in $D_{1}$ assigned to $z$ 
in \eqref{(z1def)} is a minimum $c_{1}$-cost $(f_{1},g_{1})$-bounded $m$-flow of $D_{1}$.  
\finbox 
\end{corollary}

Let $z$ be an $(f,g)$-bounded $m$-flow and let $D_{z}$ be 
the usual auxiliary digraph belonging to $z$.  
The sets of forward and backward edges in $F_{z}$ are denoted by 
$F_{{\bf f} }$ and $F_{{\bf b} }$, respectively.  
Let 
\[
L_{{\bf f} }:= \{uv\in F_{{\bf f} }:  uv\in L, z(uv) = \beta -1   \},
\quad
L_{{\bf b} }:= \{uv\in F_{{\bf b} }:  vu\in L,  z(vu) =\beta \}.
\]

\begin{lemma} \label{betaminchar} 
An integral $(f,g)$-bounded $m$-flow $z$
with $\max \{z(e):e\in L\}\leq \beta$ minimizes 
the number of $\beta$-valued (that is, $g$-saturated) elements of $L$ 
if and only if,
in every di-circuit of $D_{z}$, 
the number of $L_{\bf b}$-edges  is at most the number of $L_{\bf f}$-edges.  
\end{lemma}

\Proof
Suppose first that $z$ is an integral $(f,g)$-bounded $m$-flow
for which the auxiliary digraph $D_{z}$ belonging to $z$ includes a
di-circuit $C_{z}$ which has more $L_{\bf b}$-edges than $L_{\bf f}$-edges.  
Let $C$ denote the circuit of $D$ corresponding to $C_{z}$
(that is, $C$ is obtained from $C_{z}$ by reversing the backward edges of $C_{z}$).  
Define $z'$ as follows:
\begin{equation} 
z'(uv):= \begin{cases}
  z(uv)+1 & \hbox{if $uv\in C_{z}$ is a forward edge}, \cr 
z(uv)-1 & \hbox{if $vu\in C_{z}$ is a backward edge}, \cr
z(uv) & \hbox{if $uv\in A-C$}.  \end{cases} 
\label{(z'def)} 
\end{equation}
\noindent 
Then $z'$ is an integral $(f,g)$-bounded $m$-flow that
saturates less $L$-edges than $z$ does.

To see the converse, suppose
that $z$ is an integral $(f,g)$-bounded $m$-flow 
for which the number of $\beta$-valued (that is, saturated) $L$-edges is 
not minimum.

Consider the digraph $D_{1}$ defined above 
along with the bounding functions 
$(f_{1},g_{1})$ on its edge-set in \eqref{f1g1def}.
Let $z_{1}$ be the $(f_{1},g_{1})$-bounded $m$-flow 
assigned to $z$ in \eqref{(z1def)}.
By Lemma~\ref{minbeta}, $z_{1}$ is not a minimum $c_{1}$-cost
$(f_{1},g_{1})$-bounded $m$-flow in $D_{1}$.  
By applying Lemma~\ref{javkor} to $x:=z_{1}$, 
we obtain that the auxiliary digraph $D_{x}$ belonging to $x$ 
includes a di-circuit $C_{x}$ whose $c_{x}$-cost is negative.

Let $e=uv$ be an edge of $L$.  
Recall that, to define $D_{1}$, we added a new edge $e'$ parallel to $e$.  
Let $e''=vu$ be the edge arising from $e'$ by reversing it.  
Then we have the following equivalences:
\begin{align*}
z(e) = \beta - 1 
& \Leftrightarrow 
uv\in L_{\bf f}\subseteq A_{z}
\\ & \Leftrightarrow 
\mbox{$e'$ is a forward edge in $D_{x}$ with $c_{x}(e')=1$},
\\
z(e)=\beta
\phantom{{}-1}  
& \Leftrightarrow
 vu\in L_{\bf b}\subseteq A_{z}
\Leftrightarrow 
z_{1}(e') = 1
\\ & \Leftrightarrow 
\mbox{$e''=vu$ is a backward edge in $D_{x}$ (and hence $c_{x}(e'')=-1)$}.
\end{align*}
In addition, the $c_{x}$-cost of the edges (forward or backward) 
associated with $e$ with $z(e)< \beta -1$  
is equal to zero.
These observations imply that the negative di-circuit $C_{x}$
(with respect to $c_{x}$) in $D_{x}$ defines a di-circuit of $D_{z}$ which
contains more $L_{\bf b}$-edges than $L_{\bf f}$-edges.  
\finbox

\subsection{The characterization}
\label{SCcharthm}

Recall the cost-vector $\underline{c}$ defined in \eqref{(cdef)},
which is a $k$-dimensional vector with $k\leq 2 | F |$. 
The main result of Section~\ref{SCgall} is as follows.

\begin{theorem}   \label{gall.main} 
For an element $z\in \odotZ{Q}=\odotZ{Q}(f,g;m)$, 
the following properties are equivalent.

\noindent {\rm (A)} 
\ $z$ is decreasingly minimal on $F$.

\noindent {\rm (B)} 
\ There is no $z$-improving di-circuit in the auxiliary digraph $D_{z}$.

\noindent {\rm (C)} 
\ There is an integer-valued potential-vector function 
$\underline{\pi} $ on $V$ which is $\underline{c}$-feasible
in $D_{z}$,
that is, $\underline{\pi} (v) - \underline{\pi} (u) \preceq \underline{c}(uv)$ 
for every edge $uv \in A_{z}$,
where the dimension of $\underline{\pi} $ is bounded by $2|F|$.
\end{theorem}

\Proof
For the proof it is convenient to highlight the condition:

\medskip

{\rm (B$'$)} 
\ There is no di-circuit $C$ with $\widetilde {\underline{c}}(C)\prec \underline{0}_{k}$
in the auxiliary digraph $D_{z}$.

\medskip

\noindent
Lemma~\ref{LMekvi} shows 
the equivalence of (B) and (B$'$),
whereas the equivalence of (B$'$) and (C) is shown in Theorem~\ref{ext-Gallai}.
The implication
``(A) $\Rightarrow$ (B)''
is obvious from the definition,
and now we turn to the proof of ``(B) $\Rightarrow$ (A).''

Let $z$ be an $(f,g)$-bounded integral $m$-flow 
for which there is no $z$-improving di-circuit in the auxiliary digraph $D_{z}$.
To derive that $z$ is $F$-dec-min, we use induction on $\vert F\vert$. 
As $z$ is $F$-dec-min when $F$ is empty, we assume that 
$\vert F\vert \geq 1$.  
We can assume that $F$ contains no $(f,g)$-tight edges,
 since taking out an $(f,g)$-tight edge from $F$ affects neither
the set of $z$-improving di-circuits, nor the $F$-dec-minimality of $z$.

Let $\beta :=\max\{z(e):e\in F\}$.  
Then $\max\{z'(e):e\in F\}\leq \beta$ 
holds for any $F$-dec-min member $z'$ of $\odotZ{Q}$, 
therefore we can assume that $\beta =\max\{g(e):e\in F\}$.  
Let $L:=\{e\in F, g(e)=\beta \}$.

Since $D_{z}$ admits no $z$-improving di-circuit, 
it follows, 
in particular, that there is no di-circuit containing more 
$L_{\bf b}$-edges than $L_{\bf f}$-edges.  
By Lemma~\ref{betaminchar}, 
$z$ minimizes the number of $F$-edges 
with $z(e)=\beta$, and this means that $z$ is pre-dec-min on $F$.

Consider the chain ${\cal C}'$ used in Theorem~\ref{Fdecmin} along with the
definition of $(f',g')$ given in 
\eqref{(f'g'def1)} and \eqref{(f'g'def2)}.  
By (the proof of) Theorem~\ref{reduction}, $z$ is $(f',g')$-bounded.  
Recall that $L'$ was defined before Claim \ref{claimL'}
to be the subset of $L$ consisting of those elements of $L$ 
that enter at least one member of ${\cal C}'$, 
while we defined $F':=F-L'$.  
We pointed out that $L'$ is non-empty, that is, 
$F'$ is a proper subset of $F$.  
Furthermore the definitions of $(f',g')$ and $L'$ imply that 
every edge in $A-L$ leaving or entering a member of ${\cal C}'$ 
is $(f',g')$-tight,
every edge in $L$ leaving a member of ${\cal C}'$ is 
$(f',g')$-tight,
and  every edge in $L$ entering at least two members of ${\cal C}'$ is 
$(f',g')$-tight.

\medskip

Let $D'_{z}$ denote the auxiliary digraph belonging to $z$ with respect to $(f',g')$.  
Because $(f',g')$-tight edges of $D$ do not define any edge of $D'_{z}$, 
we conclude that, for any member $C_{i}$ of ${\cal C}'$,
if $e=uv$ is a forward edge of $D'_{z}$ entering $C_{i}$, then
$f'(e)=\beta -1$, $g'(e)=\beta$, 
and $e$ does not enter any other member of ${\cal C}'$.  
Analogously, if $e=uv$ is a backward edge of $D'_{z}$ leaving $C_{i}$, 
then $f'(vu)=\beta -1$, $g'(vu)=\beta$, and
$e=uv$ does not leave any other member of ${\cal C}'$.  
It follows for any di-circuit $K'$ of $D'_{z}$ that, 
if $K$ denotes the circuit of $D$ corresponding to $K'$, 
then the number of $F$-edges $e$ of $K$ with $z(e)=\beta -1$ 
entering $C_{i}$ is equal to the number of $F$-edges of $K$ 
with $z(e)=\beta$ leaving $C_{i}$.  
This implies that if $K'$ is a
$z$-improving di-circuit of $D'_{z}$ with respect to $F'$, 
then $K'$ is $z$-improving di-circuit in $D_{z}$ with respect to $F$.

By our hypothesis, $D_{z}$ includes no $z$-improving di-circuit, and therefore 
$D'_{z}$ includes no $z$-improving di-circuit with respect to $F'$, either.  
Since $\vert F'\vert <\vert F\vert $, we conclude by induction that 
$z$ is $F'$-dec-min with respect to $(f',g')$, 
implying, via Theorem~\ref{reduction},
that $z$ is $F$-dec-min.  
\finbox

\medskip

\begin{remark} \rm \label{RMpotveccomp2}
As we applied Theorem~\ref{ext-Gallai} 
for proving implication ``(B) $\Rightarrow $ (C)'' in Theorem~\ref{gall.main}
and, in the present case, we have $k \leq 2 | F |$ for the $k$-dimensional
cost-vector $\underline{c}$ defined in \eqref{(cdef)},
we can conclude, 
by Remark~\ref{RMpotveccomp1},
that the potential-vector $\underline{\pi}$ occurring in (C) can be
computed in strongly polynomial time
for a given $F$-dec-min element  $z \in \odotZ{Q}$.
\finbox
\end{remark}

\begin{remark} \rm  \label{char.link} 
From a theoretical computer science
point of view, a slight drawback of the characterization in Theorem~\ref{MAIN} 
is that, in order to be convinced that $z$ is indeed
$F$-dec-min, one must believe the correctness of $(f\sp{*},g\sp{*})$.
In this respect, Property {\rm (C)} in Theorem~\ref{gall.main} is more
convincing since it provides a certificate for $z$ to be $F$-dec-min
whose validity can be checked immediately.

Just for an analogy to understand better this aspect of certificates,
consider the well-known maximum weight perfect matching problem in a
bipartite graph $G=(S,T;E)$ endowed with a weight-function $w$ on $E$.
On one hand, one can prove the characterization that there is a
subgraph $G\sp{*}=(S,T;E\sp{*})$ 
of $G$ such that a perfect matching $M$ of $G$
is of maximum $w$-weight if and only if 
$M\subseteq E\sp{*}$.  
(This result intuitively corresponds to Theorem~\ref{MAIN}).  
This certificate $E\sp{*}$, however, is convincing
 (for the optimality of $M$) only if we
can check that it has been correctly computed.  
On the other hand,
Egerv\'ary's classic theorem provides an immediately checkable
certificate for $M$ to be of maximum $w$-weight:  a function 
$\pi: S\cup T \rightarrow {\bf R} $ 
for which $\pi (s)+\pi (t)\geq w(st)$
for every edge $st\in E$ and $\pi (s)+\pi (t)=w(st)$ for every edge $st\in M$. 
(This result intuitively corresponds to the equivalence of
{\rm (A)} and {\rm (C)} in Theorem~\ref{gall.main}).  
\finbox
\end{remark}

\section{Existence of an $F$-dec-min $m$-flow} 
\label{veges}

In the previous sections, we assumed that the bounding functions $f$
and $g$ were finite-valued on $F$.  
In the more general case, where we
allow edges in $F$ as well to have $f(e)=-\infty $ or $g(e)=+\infty$,
it may occur that no 
$F$-dec-min feasible $m$-flow exists at all. 
 For example, if $D$ is a di-circuit, 
$F=A$, $m\equiv 0$, $f\equiv -\infty$, and $g\equiv 0$, 
then $z\equiv k$ is a feasible $m$-flow for each integer $k\leq 0$, 
implying that in this case there is no $F$-dec-min feasible $m$-flow.  
The main goal of this section is to describe a characterization 
for the existence of an $F$-dec-min feasible $m$-flow.  
As a consequence of this characterization, we show how
Theorem~\ref{MAIN} and its algorithmic approach can be extended to
this more general case.

As before, let $D=(V,A)$ be a digraph and $F\subseteq A$ a non-empty subset of edges.  
Let $m:V \rightarrow {\bf Z} $ be a function on $V$
and let $f:A \rightarrow {\bf Z}\cup \{-\infty \}$ and 
$g:A \rightarrow {\bf Z}\cup \{+\infty \}$ 
be bounding functions on $A$ such that there
is a feasible (that is, $(f,g)$-bounded) $m$-flow in $D$.  
Recall that
$\odotZ{Q}(f,g;m)$ denoted the set of integral $(f,g)$-bounded $m$-flows.  
In what follows, all the occurring functions (bounds, flows) 
are assumed to be integer-valued even if this is not mentioned explicitly.

We start by exhibiting an easy reduction by which we can assume that
$g$ is finite-valued on $F$.

\begin{lemma} \label{gfinonF} 
There is a function $g'$ on $A$ which is
finite-valued on $F$ such that the (possibly empty) set of $F$-dec-min
elements of $\odotZ{Q}:=\odotZ{Q}(f,g;m)$ is equal to the set of
$F$-dec-min elements of $\odotZ{Q'}:=\odotZ{Q}(f,g';m)$. 
\end{lemma}

\Proof 
Let $z_{1}$ be an element of $\odotZ{Q}$ and let $\beta$ denote
the maximum value of its components in $F$.
Define $g'$ as follows:
\begin{equation} 
g'(e):= \begin{cases} 
 \min \{g(e),\beta \} & \ \ \hbox{if}\ \ \ e\in F, \ 
  \cr
  g(e) & \ \ \hbox{if}\ \ \ e \in A-F. 
        \end{cases}
\label{(g'def)} 
\end{equation}
\noindent
As $g'\leq g$, we have $\odotZ{Q'} \subseteq \odotZ{Q}$.
In particular, an $F$-dec-min element $z'$ of $\odotZ{Q'}$ 
is in $\odotZ{Q}$, and we claim that $z'$ is actually $F$-dec-min in $\odotZ{Q}$.  
Indeed, if we had an element $z'' \in \odotZ{Q}$ which is decreasingly smaller on $F$
than $z'$, then $z''$ is not in $\odotZ{Q'}$, that is, $z''$ is not $(f,g')$-bounded.  
Therefore there is an edge $a\in F$ for which $z''(a)>\beta$, 
implying that 
$\max \{z''(e) :  e\in F\} >\beta = \max \{z'(e):e\in F\}$.  
But this contradicts the assumption that
$z''$ is decreasingly smaller on $F$ than $z'$.

Conversely, suppose that $z$ is an $F$-dec-min element of $\odotZ{Q}$.
Since the largest component of $z_{1}$ in $F$ is $\beta$, 
the largest component of $z$ in $F$
is at most $\beta$, and hence $z\in \odotZ{Q'}$.
This and $\odotZ{Q'} \subseteq \odotZ{Q}$ imply that $z$ is an
$F$-dec-min element of $\odotZ{Q'}$.  
\finbox \medskip

\begin{theorem}   \label{findecmin} 
Let $D=(V,A)$ be a digraph and $F\subseteq A$ a subset of edges. 
Let $m:V \rightarrow {\bf Z} $ be a function on $V$ and let 
$f:A \rightarrow {\bf Z}\cup \{-\infty \}$ and
$g:A \rightarrow {\bf Z}\cup \{+\infty \}$ be bounding functions on $A$
such that there is a feasible (that is, $(f,g)$-bounded) $m$-flow in $D$.  
Define digraph $D\sp{\infty }=(V,A\sp{\infty})$ by
\begin{equation} 
A\sp{\infty} := \{e:  e\in A, f(e)=-\infty \} 
 \cup \{vu:  uv\in A-F, g(uv)=+\infty \}.  
\label{(Ainfdef)} 
\end{equation} 
The following properties are equivalent.

\noindent 
{\rm (A)} \ 
There exists an $F$-dec-min $(f,g)$-bounded integral $m$-flow.

\noindent 
{\rm (B)} \ 
There is no di-circuit $C$ in $D\sp{\infty}$ with $C \cap F \ne \emptyset$.

\noindent 
{\rm (C)} \ 
Each edge $e\in F$ with $f(e)=-\infty $ enters a
subset $S_{e}$ for which $\delta_{A\sp{\infty}}(S_{e}) = 0$. 
\end{theorem}  

\Proof
Since each of the three properties holds when $F=\emptyset$, 
we can assume that $F$ is non-empty.  
As a first step, we make the upper bound function $g$ finite-valued on $F$.

\begin{claim} \label{gfinonF.1} 
The theorem follows from its special case
when $g(e)$ is finite for each $e\in F$. 
\end{claim}  

\Proof 
Consider the function $g'$ introduced in \eqref{(g'def)}, 
and suppose that the theorem holds when $g$ is replaced by $g'$.  
To derive the theorem for the original $g$, observe first that changing
$g$ to $g'$ does not affect the digraph $D\sp{\infty}$ 
because $g'$ may differ from $g$ only on the elements of $F$.  
Since both Property (B) and Property (C) depend only on $D\sp{\infty}$, 
these properties are not affected by replacing $g$ with $g'$, and hence they are
equivalent (with respect to $g$).  
Furthermore, Lemma~\ref{gfinonF}
implies that Property (A) holds with respect to $g$ precisely 
if it holds with respect to $g'$.  
\finbox \medskip

By the claim, we can assume that $g$ is finite-valued on $F$. 
Note that in this case
\begin{equation}
 A\sp{\infty} = \{e:  e\in A, f(e)=-\infty \} \cup \{vu:  uv\in A, g(uv)=+\infty \}. 
 \label{(Ainfdef.1)} 
\end{equation}

(A) $\Rightarrow $ (B) \ Let $z$ be an $F$-dec-min element of 
$\odotZ{Q}$.  
Suppose indirectly that $D\sp{\infty} $ includes a di-circuit $C$
intersecting $F$.  
For $uv\in A$, define $z'(uv)$ as follows:
\begin{equation} 
\label{existdecminprf1}
z'(uv):= \begin{cases} 
z(uv)- 1 
& \ \ \hbox{if}\ \ \ uv\in C, \ uv\in A ,
 \\
z(uv)+ 1 
& \ \ \hbox{if}\ \ \ vu\in C, \ vu \in A-F, 
 \\
z(uv) & \ \ \hbox{otherwise}. 
   \end{cases} 
\end{equation}
\noindent
Then $z'$ is also a feasible $m$-flow in $D$, which is
decreasingly smaller on $F$ than $z$, a contradiction.

\medskip

(B) $\Rightarrow $ (C) \ For any edge $e=ts\in F$ with $f(e)=-\infty$, let 
$S_{e}$ denote the set of nodes which are reachable from $s$ in $D\sp{\infty}$. 
Then $e$ enters $S_{e}$ since if we had 
$t \in S_{e}$, 
then there is an $st$-dipath $P$ in $D\sp{\infty}$, 
and the di-circuit $C=P+e$ would violate Property (B).

\medskip

(C) $\Rightarrow $ (A) First we provide a condition for an edge $e\in
F$ which ensures that $z(e)$ cannot be arbitrarily small for $z\in \odotZ{Q}$.

\begin{claim} \label{l-bound.e0} 
Let $S\subset V$ be a set for which $\delta_{A\sp{\infty} }(S)=0$, 
and let $e_{0}\in F$ entering $S$.   
Then, for any $(f,g)$-feasible $m$-flow $z$,
\begin{equation} 
z(e_{0}) \geq \widetilde m(S) - [ \varrho_{g}(S) - g(e_{0}) ] +  \delta_{f}(S), 
\label{(zlbound)} 
\end{equation} 
and the right-hand side is finite.  
\end{claim}

\Proof 
Since $z\leq g$ and $e_{0}$ enters $S$, we have 
\[
\varrho_{z}(S) - z(e_{0}) \leq \varrho_{g}(S) - g(e_{0}) ,
\] 
from which
\[
  \widetilde m(S)= \varrho_{z}(S)- \delta_{z}(S)
 = z(e_{0}) + [\varrho_{z}(S) - z(e_{0})] - \delta_{z}(S) \leq z(e_{0}) 
  + [ \varrho_{g}(S) - g(e_{0})] - \delta_{f}(S),
\] 
implying \eqref{(zlbound)}.

Furthermore, $\delta_{A\sp{\infty} }(S)=0$ implies that $f(e)>-\infty$ 
for every edge $e$ of $D$ leaving $S$ and that $g(e) <+\infty $ 
for every edge $e$ of $D$ entering $S$, from which the finiteness of 
the right-hand side of \eqref{(zlbound)} follows.  
\finbox \medskip

Assume indirectly that no $F$-dec-min $(f,g)$-bounded $m$-flow exists,
that is, for every $(f,g)$-bounded $m$-flow, there exists another one
which is decreasingly smaller on $F$.  
This implies that there is an
edge $e_{0}$ in $F$ for which there is an $(f,g)$-bounded $m$-flow
with $z(e_{0})\leq K$ for an arbitrarily small integer $K$.
By Claim \ref{l-bound.e0}, 
$e_{0}$ cannot enter any subset $S\subset V$ with 
$\delta_{A\sp{\infty}}(S)=0$, 
contradicting Property (C).  
This completes the proof of Theorem~\ref{findecmin}.
\finbox \finboxHere

\begin{corollary} \label{fgveges} 
Let $Q=Q(f,g;m)$ be the set of $(f,g)$-bounded $m$-flows.  
If $\odotZ{Q}$ has an $F$-dec-min element, 
then there are bounding functions $(f',g')$  for which 
the sets of $F$-dec-min elements of $\odotZ{Q}(f',g';m)$ and of $\odotZ{Q}$ are the same, 
and both $f'$ and $g'$ are finite-valued on $F$.  
\end{corollary}

\Proof 
Lemma~\ref{gfinonF} implies that the upper-bound function $g'$
defined in \eqref{(g'def)} is finite-valued on $F$, and replacing $g$
by $g'$ does not affect the set of $F$-dec-min elements.  
Since $\odotZ{Q}$ has an $F$-dec-min element, 
Theorem~\ref{findecmin} implies that every edge $e\in F$ 
with $f(e)=-\infty $ enters a subset $S_{e}$
for which $\delta _{A\sp{\infty} }(S_{e}) = 0$.  
This and Claim \ref{l-bound.e0} imply, that there is a finite lower bound 
\begin{equation}  \label{(f'(e))} 
f'(e):= \widetilde m(S_{e}) - [ \varrho_{g}(S_{e}) - g(e) ] + \delta_f(S_{e}) .
\end{equation} 
For these $f'$ and $g'$, the set of $F$-dec-min
elements of $\odotZ{Q}$ is the same as the set of $F$-dec-min elements
of $\odotZ{Q}(f',g';m)$.  \finbox

\medskip 

Corollary~\ref{fgveges} implies that our main theorem
(Theorem~\ref{MAIN}) holds almost word for word in the general
case when $(f,g)$ is not assumed to be finite-valued on $F$:  
the only difference is that the existence of an $F$-dec-min element of $\odotZ{Q}$ 
must be assumed.

\begin{theorem}   \label{MAINb} 
Let $D, F, f, g, m$ be the same as in Theorem~\ref{findecmin}, 
and let $Q=Q(f,g;m)$ be the set of $(f,g)$-bounded feasible $m$-flows.  
Assume that there exists an $F$-dec-min element of $\odotZ{Q}$.  
Then there exists a pair $(f\sp{*},g\sp{*})$ of
integer-valued functions on
$A$ with $f\leq f\sp{*} \leq g\sp{*} \leq g$ 
(allowing $f\sp{*}(e)=-\infty $ and $g\sp{*}(e)=+\infty $ for $e\in A-F$)
such that an integral $(f,g)$-bounded $m$-flow $z$ 
is decreasingly minimal on $F$ 
if and only if 
$z$ is an integral $(f\sp{*},g\sp{*})$-bounded $m$-flow.
Moreover, the box $T(f\sp{*},g\sp{*})$ is narrow on $F$ in the sense
that $0\leq g\sp{*}(e)-f\sp{*}(e)\leq 1$ for every $e\in F$.  
\finbox
\end{theorem}

We mention that the description above immediately gives rise to a
strongly polynomial algorithm that terminates by providing either a
di-circuit $C$ in $D\sp{\infty}$ intersecting $F$ 
(which is a certificate for the non-existence of an $F$-dec-min element) 
or else the bounding functions $(f',g')$ 
occurring in Corollary~\ref{fgveges}
which are finite-valued on $F$.  
The only subroutine needed here is the one
to compute the set $S_{e}$ of nodes reachable in $D\sp{\infty}$
from a specified node.  
This can easily be realized, for example, by a breadth-first search.

\section{Computing an $L$-upper-minimizer $m$-flow and the dual optimal chain} 
\label{algo5}

In the previous sections we provided a necessary and sufficient condition 
for the existence of an $F$-dec-min integral $(f,g)$-bounded $m$-flow, 
characterized these $m$-flows, and described their set 
as the set of integral $(f\sp{*},g\sp{*})$-feasible $m$-flows.
Our next goal is to consider algorithmic questions and construct
strongly polynomial algorithms for the results developed earlier.

In the present section, we describe an alternative, algorithmic proof of 
Theorems \ref{minL} and \ref{optkritx}.  
This algorithm will actually be used in the special case, 
described in Theorem \ref{Fdecmin}, 
for computing the dual optimal chain ${\cal C'}$ 
characterizing the $F$-pre-dec-min elements of $\odotZ{Q} = \odotZ{Q}(f,g;m)$.  
This chain, as described in Theorem~\ref{reduction}, 
immediately gives rise to a tightening $(f',g')$ of $(f,g)$ and
a proper subset $F'$ of $F$ 
with the property that the set of $F$-dec-min elements 
of $\odotZ{Q}$ is the same as 
the set of $F'$-dec-min elements of $\odotZ{Q'} = \odotZ{Q}(f',g';m)$.

In the light of this algorithmic proof,
the original proof of Theorems \ref{minL} and \ref{optkritx}
may seem superfluous, but we keep both proofs
because the one in 
Section~\ref{flow3} is more transparent and technically simpler 
than the algorithmic approach to be presented below.

The algorithm computes an integer-valued $L$-upper-minimizer
$(f,g)$-bounded $m$-flow as well as a maximizer chain $\cal C$ in
\eqref{(minmaxL)} meeting the optimality criteria in Theorem~\ref{optkritx}.  
As before, $D=(V,A)$ is a digraph and we assume that
$L$ is a subset of $A$ for which $-\infty <f(e)<g(e)<+\infty $ 
for each edge $e\in L$.  
(For edges in $A-L$, $f(e)=-\infty $ and $g(e)=+\infty$ are allowed.)  
Our primal goal is to find an integral
$(f,g)$-bounded $m$-flow $g$-saturating a minimum number of elements of $L$.  
To this end, 
we apply the technique used already in Section~\ref{SCminsatedge}
which relies on cheapest feasible flows.  
However, these two frameworks differ in the following respects.
In Section~\ref{SCminsatedge}, $\beta $ and $F$ played a role,
while these parameters do not occur here.  
Another difference is that in Section~\ref{SCminsatedge} 
we relied only on the primal optimum of the min-cost flow problem, 
while here it is central to compute the dual optimum, as well.

Similarly to the approach of Section~\ref{SCminsatedge},
we introduce a parallel copy $e'$ for each element $e\in L$.  
Let $L'$ denote the set of new edges.  
We shall refer to the edges in $A$ as old or original edges.  
Let $A_{1}:=A\cup L'$, $D'=(V,L')$, and $D_{1}=(V,A\cup L')$.  
Define $g\sp{-}$ on $A$ by $g\sp{-}:=g-\chi_{L}$, 
that is, we reduce $g(e)$ by 1 for each $e\in L$.
Let $f_{1}$ and $g_{1}$ be bounding functions on $A_{1}$ defined by \eqref{f1g1def},
and $c_{1}$ be a $(0,1)$-valued cost-function on $A_{1}$ defined by \eqref{c1def}.

Our goal is to find an $(f,g)$-bounded integer-valued $m$-flow in $D$
admitting a minimum number of $g$-saturated $L$-edges.  
We claim that this problem is equivalent to finding a minimum $c_{1}$-cost
$(f_{1},g_{1})$-bounded integer-valued $m$-flow in $D_{1}$.  
Indeed, let $z$ be an $(f,g)$-bounded $m$-flow in $D$ and let $X:=\{e\in L:
z(e)=g(e)\}$ be the set of $g$-saturated members of $L$.  
Let $X'$ denote the subset of $L'$ corresponding to $X$.  
Define an $m$-flow $z_{1}$ in $D_{1}$ as follows:
\[
 z_{1}(e):= \begin{cases}
 z(e) & \hbox{if $e\in A-X$}, 
  \cr 
  g(e)-1 & \hbox{if $e\in X$}, 
  \cr
  1 & \hbox{if $e\in X'$}, 
  \cr 
  0 & \hbox{if $e\in L'-X'$.  }\ 
         \end{cases}
\]
\noindent
Then $z_{1}$ is an $(f_{1},g_{1})$-bounded $m$-flow in $D_{1}$ whose
$c_{1}$-cost is $\vert X\vert $. 
Conversely, let $z_{1}$ be a minimum cost integer-valued $(f_{1},g_{1})$-bounded $m$-flow in $D_{1}$.  
Observe that if
$z_{1}(e')=1$ for some $e'\in L'$, then $z_{1}(e)=g_{1}(e)=g(e)-1$ where $e$
is the edge in $L$ corresponding to $e'$.  Indeed, if we had
$z_{1}(e)\leq g(e)-2$, then the $m$-flow obtained from $z_{1}$ by adding 1
to $z_{1}(e)$ and subtracting 1 from $z_{1}(e')$ would be of smaller cost.
It follows that the $m$-flow $z$ in $D$ defined by
\begin{equation} 
z(e):= \begin{cases} 
 z_{1}(e)+z_{1}(e') & \hbox{if $e\in L$}, 
\cr
 z_{1}(e) & \hbox{if $e\in A-L$ }\ 
   \end{cases} 
\label{(xdef)}
 \end{equation}
is an $(f,g)$-bounded $m$-flow in $D$, for which the number
of $g$-saturated $L$-edges is exactly the $c_{1}$-cost of $z_{1}$.

Therefore, we concentrate on finding an integer-valued min-cost
$(f_{1},g_{1})$-bounded $m$-flow in $D_{1}$.  In order to describe the dual
optimization problem, 
let $N$ denote the node-edge signed incidence matrix of $D$, 
that is, the entry of $N$ corresponding to a node $u$
and to an edge $e\in A$ is 1 if $e$ enters $u$, 
 $-1$ if $e$ leaves $u$, and 0 otherwise.  
Let $N'$ denote the analogous signed incidence
matrix of $D'$, and let $N_{1}=[N,N']$.  
Note that $N_{1}$ is the signed
incidence matrix of $D_{1}$ and hence it is totally unimodular.  
The primal linear program is as follows:
\begin{equation} 
 \min \{c_{1}z_{1} :  \ N_{1}z_{1}=m, \ z_{1}\geq f_{1}, \ -z_{1} \geq -g_{1}\}.
 \label{(primalflow)} 
\end{equation}

The dual linear program is as follows:
\begin{equation} 
\max \{ym + v_{1}f_{1} - w_{1}g_{1}: 
 \ yN_{1} + v_{1} -w_{1} = c_{1}, \ v_{1}\geq 0, \ w_{1}\geq 0\}.  
\label{(dualflow)} 
\end{equation}
\noindent
Note that the components of $v_{1}=(v,v')$ correspond to the
edges in $A$ and in $L'$, respectively, and the analogous statement
holds for $w_{1}=(w,w').$ Since $N_{1}$ is totally unimodular, both the
primal and the dual optimal solution can be chosen integer-valued.

If $(y,v_{1},w_{1})$ is a dual solution and both $v_{1}(e)$ and $w_{1}(e)$ are
positive on an edge $e\in A_{1}$, 
then reducing both $v_{1}(e)$ and $w_{1}(e)$ by 
their minimum $\delta := \min \{v_{1}(e), w_{1}(e)\}$,
we obtain another dual solution
whose dual cost is larger by 
$\delta (g_{1}(e)-f_{1}(e)) \geq 0$
than the dual cost
$ym + v_{1}f_{1} - w_{1}g_{1}$ of $(y,v_{1},w_{1})$.  
Therefore it suffices to
consider only those optimal dual solutions $(y,v_{1},w_{1})$ for which
$\min \{v_{1}(e), w_{1}(e)\}=0$ for every edge $e\in A_{1}$.  Observe that
for such an optimal dual solution $(y,v_{1},w_{1})$, 
since $v_{1}$ and $w_{1}$ are non-negative, $y$ 
uniquely determines $v_{1}$ and $w_{1}$.  
Namely,
for an edge $e=s t\in A$, we have $c_{1}(e)=0$ and hence
\begin{align} 
v_{1}(e)& := \begin{cases}
   0 & \hbox{if $y(t)-y(s)\geq 0$}, 
   \cr y(s)- y(t) & \hbox{if $y(t)-y(s)<0$, }\ 
        \end{cases} 
\label{(z1A)} 
\\
w_{1}(e) &:= \begin{cases}
0 & \hbox{if $y(t)-y(s)\leq 0$}, 
 \cr
 y(t)- y(s) & \hbox{if $y(t)-y(s)>0$}. 
        \end{cases} 
\label{(w1A)} 
\end{align} 
\noindent
For an edge $e'=s t\in L'$, we have $c_{1}(e')=1$ and hence
\begin{align} 
 v_{1}(e') & := \begin{cases}
   0 & \hbox{if $y(t)-y(s)\geq 1$},
  \cr 
 y(s)- y(t)+ 1 & \hbox{if $y(t)-y(s)<1$}, 
              \end{cases} 
\label{(z1L)}
\\
w_{1}(e') &:= \begin{cases}
   0 & \hbox{if $y(t)-y(s) \leq 1$}, 
   \cr
 y(t)- y(s) -1 & \hbox{if $y(t)-y(s)>1.$}\ 
           \end{cases}
 \label{(w1L)}
\end{align}

Let $z_{1}$ be an integer-valued primal optimum, that is, $z_{1}$ is a
minimum $c_{1}$-cost $(f_{1},g_{1})$-bounded $m$-flow in $D_{1}$.  
Let $z$ be the $(f,g)$-bounded $m$-flow in $D$ defined in \eqref{(xdef)}.  
As noted above, $z$ is $L$-upper-minimizer.  
Let $(y,v_{1},w_{1})$ be an integer-valued dual optimum.

Note that the minimum cost flow algorithm of 
Ford and Fulkerson \cite{Ford-Fulkerson} 
computes a minimum-cost feasible flow of given
amount along with the optimal dual solution.  
This algorithm relies on a max-flow algorithm as a subroutine.  
If one uses the strongly polynomial max-flow algorithm 
of Edmonds and Karp \cite{Edmonds-Karp},
that is, if the augmentation is made always along a shortest path in
the corresponding auxiliary digraph, and, furthermore, if the
cost-function is $(0,1)$-valued, then the min-cost flow algorithm of
Ford and Fulkerson is strongly polynomial. 
 (In other words, we do not
need to use a more sophisticated strongly polynomial algorithm---the
first one found by Tardos \cite{Tardos1}---for the general min-cost
flow problem when the cost-function is arbitrary.)  
With a standard reduction technique, 
the min-cost flow algorithm of Ford and Fulkerson can easily 
be transformed to one for computing a feasible min-cost $m$-flow.  
Therefore, we conclude that the integer-valued optimal
solutions to the primal and dual linear programs above can be computed
in strongly polynomial time via the Ford-Fulkerson min-cost flow algorithm.

Since $\widetilde m(V)=0$, by adding a constant to the components of
$y$, we obtain another optimal dual solution.  Therefore we may assume
that the smallest component of $y$ is 0. 
Let $0=y_{0}<y_{1}<y_{2}<\cdots <y_{q}$ 
be the distinct values of the components of $y$, and consider
the chain of subsets 
$V_{1}\supset V_{2}\supset \cdots \supset V_{q}$ of $V$
where $V_{i}:=\{u \in V:  y(u)\geq y_{i}\}$. 
 (In the special case when $y\equiv 0$, the chain in question is empty, that is, $q=0$).

Note that 
\begin{equation} 
 ym = \sum_{i=1}\sp{q} (y_{i}-y_{i-1}) \widetilde m(V_{i}).
\label{(atirt)} 
\end{equation}
\noindent
We may assume that the difference of subsequent $y_{i}$ values is 1. 
Indeed, if $y_{i+1}-y_{i}\geq 2$ for some $i$, then by subtracting 1 
from $y(u)$ for each $u \in V_{i+1}$, by subtracting 1 from $v_{1}(e)$
for each $e\in A_{1}$ leaving $V_{i+1}$, and by subtracting 1 from
$w_{1}(e)$ for each $e\in A_{1}$ entering $V_{i+1}$, 
we obtain another dual feasible solution $(y',v_{1}',w_{1}')$.  
By \eqref{(atirt)}, $y'm = ym - \widetilde m(V_{i+1})$.  
For the revised $v_{1}'$ and $w_{1}'$, we have
\begin{align*} 
v_{1}' f_{1} &= v_{1}f_{1} - \delta_{f_{1}}(V_{i+1}) = v_{1}f_{1} - \delta_{f}(V_{i+1}),
\\
w'_{1}g_{1} &= w_{1}g_{1} -\varrho_{g_{1}}(V_{i+1}) = w_{1}g_{1} - \varrho_{g}(V_{i+1}).
\end{align*} 
\noindent
Therefore 
\[
 y'm + v_{1}'f_{1} - w_{1}'g_{1} 
 = ym + v_{1}f_{1} - w_{1}g_{1} 
  - [\widetilde m(V_{i+1}) + \delta_{f}(V_{i+1}) - \varrho_g(V_{i+1})].
\]
\noindent
Since $\varrho_{g}(V_{i+1})-\delta_{f}(V_{i+1})\geq \widetilde
m(V_{i+1})$ by \eqref{(Hoffman)} and since $(y,v_{1},w_{1})$ is an optimal
dual solution, we obtain
\begin{align*} 
& ym + v_{1}f_{1} - w_{1}g_{1} \geq y'm + v_{1}'f_{1} - w_{1}'g_{1} 
\\ &
= ym + v_{1}f_{1}
- w_{1}g_{1} - [\widetilde m(V_{i+1}) + \delta_{f}(V_{i+1}) - \varrho
_g(V_{i+1})] \geq ym + v_{1}f_{1} - w_{1}g_{1}.
\end{align*} 
Therefore, equality must hold everywhere and hence
$(y',v_{1}',w'_{1})$ is another optimal dual solution.  
This reduction technique shows that we can assume that 
\begin{equation} 
   \hbox{ $y_{i}=i$ \ for \   $i=1,\dots ,q$}. 
\label{(yVi)} 
\end{equation}
\noindent
 Note that from an algorithmic
point of view, we get immediately the optimal dual $y$ given in \eqref{(yVi)} 
once the chain $V_{1}\supset V_{2}\supset \cdots \supset V_{q}$
belonging to an arbitrary optimal dual solution is available.

By \eqref{(yVi)}, \eqref{(z1A)}, and \eqref{(w1A)}, we have for an edge $e \in A$,
\begin{align} 
 v_{1}(e) &= \hbox{ the number of $V_{i}$'s left by $e$, }\
\label{(z1Ab)}
\\
 w_{1}(e) &= \hbox{ the number of $V_{i}$'s entered by $e$.  }\
\label{(w1Ab)} 
\end{align} 
For an edge $e'\in L'$, by \eqref{(z1L)} and \eqref{(w1L)}, we have
\begin{align} 
v_{1}(e')& = \begin{cases}0 & \hbox{if $e'$ enters a $V_{i}$}, 
 \cr
\hbox{$[$the number of $V_{i}$'s left by $e'] + 1$}\ & \hbox{if $e'$ enters no $V_{i}$, }\ 
           \end{cases} 
\label{(z1Lb)} 
\\
 w_{1}(e')& = \begin{cases}0 & \hbox{if $e'$ enters no $V_{i}$}, 
 \cr
\hbox{$[$the number of $V_{i}$'s entered by $e'] -1$}\ & \hbox{if $e'$ enters a $V_{i}$}. 
             \end{cases} 
\label{(w1Lb)}
\end{align}

The optimality criteria (complementary slackness conditions) for the
primal and dual linear programs \eqref{(primalflow)} and \eqref{(dualflow)} 
are as follows:
\begin{align} 
& \hbox{if $v_{1}(e)>0$ for some $e\in A_{1}$, then $z_{1}(e)=f_{1}(e)$},
\label{(z1optkrit)}
\\
& \hbox{if $w_{1}(e)>0$ for some $e\in A_{1}$, then $z_{1}(e)=g_{1}(e)$}.
\label{(w1optkrit)} 
\end{align}

\begin{lemma} \label{LM81crit5}
The chain $V_{1}\supset V_{2}\supset \cdots \supset V_{q}$ and the
$m$-flow $z$ defined in \eqref{(xdef)} 
meet the five optimality criteria in Theorem {\rm \ref{optkritx}}.  
Furthermore, $\varrho_{g}(V_{i})-\delta_{f}(V_{i})<+\infty $ holds 
for each $i=1,\dots ,q$.  
\end{lemma}

\Proof
 {\rm (O1)} \ Let $e\in A$ be an edge leaving a $V_{i}$.  
Then $v_{1}(e) > 0$ by \eqref{(z1Ab)}.  
By \eqref{(z1optkrit)},
$z_{1}(e)=f_{1}(e)=f(e)$,
from which $z(e)=z_{1}(e)=f(e)$ 
follows whenever $e\in A-L$.
If $e\in L$, then \eqref{(z1Lb)} implies $v_{1}(e')>0$ for the
corresponding parallel edge $e'$ in $L'$.  
By \eqref{(z1optkrit)},
$z_{1}(e')=f_{1}(e')=0$, and hence $z(e) = z_{1}(e) + z_{1}(e') = f(e)$, as
required for Criterion {\rm (O1)}.

{\rm (O2)} \ Let $e=A-L$ be an edge entering a $V_{i}$.  
Then $w_{1}(e) > 0$ by \eqref{(w1Ab)}.  
By \eqref{(w1optkrit)}, we have $z(e)= z_{1}(e) = g_{1}(e) = g(e)$, 
as required for Criterion {\rm (O2)}.

{\rm (O3)} \ Let $e\in L$ be an edge entering $V_{i}$ and let $e'$ be
the corresponding parallel edge in $L'$.  Then $w_{1}(e)>0$ by \eqref{(w1Ab)}.  
By \eqref{(w1optkrit)}, we have $z_{1}(e) = g_{1}(e)=g(e)-1$.  
Since $0=f_{1}(e') \leq z_{1}(e') \leq g_{1}(e')=1$ and $z(e)=z_{1}(e)+z_{1}(e')$, 
we obtain that $g(e)-1\leq z(e)\leq g(e)$, as
required for Criterion {\rm (O3)}.

{\rm (O4)} \ Let $e\in L$ be an edge entering at least two $V_{i}$'s,
and let $e'$ be the corresponding parallel edge in $L'$.  
By \eqref{(w1Ab)}, we have $w_{1}(e)>0$,
from which \eqref{(w1optkrit)}
implies that $z_{1}(e)=g_{1}(e)=g(e)-1$.  
By \eqref{(w1Lb)}, we have $w_{1}(e')>0$,
 from which \eqref{(w1optkrit)} implies $z_{1}(e')=g_{1}(e')=1$.
Therefore $z(e) = z_{1}(e) + z_{1}(e')=g(e)$, as required for Criterion {\rm (O4)}.

{\rm (O5)} \ Let $e\in L$ be an edge neither entering nor leaving any $V_{i}$, 
and let $e'$ be the corresponding parallel edge in $L'$.  
Since $z$ is $(f,g)$-bounded, we have $f(e)\leq z(e)$.  
By \eqref{(z1Lb)}, $v_{1}(e')=1$,
from which \eqref{(z1optkrit)} 
implies that $z_{1}(e')=f_{1}(e')=0$.  
Hence $z(e) = z_{1}(e)+z_{1}(e')\leq g_{1}(e) = g(e)-1$, 
as required for Criterion {\rm (O5)}.

To see the second part of the lemma, observe that Criterion {\rm (O1)}
implies that $\delta_{f}(V_{i})=\delta_{z}(V_{i}) >-\infty $. 
As $g(e)<+\infty $ for every edge $e\in L$, and, by Criterion {\rm (O2)}
$g(e)=v(e)<+\infty $ 
for every edge $e\in A-L$ entering $V_{i}$, 
we conclude that $\varrho_{g}(V_{i})<+\infty $, from which 
$\varrho_{g}(V_{i})-\delta_{f}(V_{i})<+\infty $, as required.  
\finbox

\medskip

Lemma~\ref{LM81crit5} and Theorem~\ref{optkritx} imply that the chain 
$V_{1} \supset V_{2} \supset \cdots \supset V_{q}$ 
computed by the algorithm above in the special case described in Theorem~\ref{reduction} 
is a dual optimal chain.
Also, the algorithm computes a minimum $c_{1}$-cost integral
$(f_1,g_1)$-bounded $m$-flow in $D_{1}$ in \eqref{(xdef)} 
and an integral $L$-upper-minimizer $m$-flow $z$ in the original digraph $D$.

\medskip

Finally, we remark that the algorithm described above can be applied
to the special case treated by Theorems \ref{Fdecmin} and \ref{reduction} 
only if the value $\beta =\beta_{F}$ defined in 
\eqref{(betaF)} is already available.  
In the next section, we show how $\beta$ can be computed efficiently.

\section{Algorithm for minimizing the largest $m$-flow value on $F$}
\label{algo1}

Our remaining algorithmic task is to describe a strongly polynomial algorithm 
for computing the smallest integer $\beta =\beta_{F}$ 
for which $\odotZ{Q}$ has an element $z$ satisfying $z(e)\leq \beta$ 
for every edge $e\in F$.  
The main tool for this computation is
the following variant of the Newton--Dinkelbach algorithm.

\subsection{Maximizing $\lceil {p(X) / b(X)}\rceil$ with 
a variant of the Newton--Dinkelbach algorithm} 
\label{SCpbmaxND}

Let $S$ be a finite ground-set.
In this section we describe a variant
of the Newton--Dinkelbach (ND) algorithm to compute the maximum 
$\lceil {p(X) / b(X)}\rceil $ 
over the subsets $X$ of $S$ with $b(X)>0$,
provided this maximum is non-negative.  
We assume that $p$ and $b$ are integer-valued set-functions 
on $S$ with $n\geq 1$ elements,
$p(\emptyset )=0$, $p(S)$ is finite
($p(X)$ may be $-\infty$ for some $X$ but it is never $+\infty $),
and $b$ is finite-valued and non-negative. 
We emphasize that there is no sign constraint on $p$ whereas $b$ is assumed to be non-negative.
The present algorithm generalizes the one described in \cite{FM21partB}
for the special case of $b(X)=\vert X\vert $, 
where $S$ is used to denote the ground-set.
In this paper, however, the algorithm will be applied to $S:=V$.

An excellent overview by Radzik \cite{Radzik} analyses 
several versions and applications of the ND-algorithm,
while a work by Goemans et al.~\cite{GGJ17} describes the most recent developments.  
We present a variant of the ND-algorithm whose specific feature is that it works
throughout with integers $\lceil {p(X) / b(X)}\rceil$. 
This has the advantage that the proof is simpler than the original one 
working with the fractions ${p(X) / b(X)}$.

The algorithm works if a subroutine is available to 
\begin{equation} 
\hbox{ find a subset of $S$ maximizing
$p(X) - \mu b(X)$ \ $(X\subseteq S)$ \ for any fixed integer $\mu \geq 0$. }\ 
\label{(ND.routine)} 
\end{equation}
\noindent
This routine will actually be needed only for special values of $\mu $ when 
$\mu =\lceil p(X)/\ell\rceil$ $\geq 0$
with $X\subseteq S$ and $1\leq \ell\leq M$, 
where $M$ denotes the largest value of $b$.  
Note that we do not have to assume that $p$ is supermodular and $b$ is submodular, 
the only requirement for the ND-algorithm is that Subroutine \eqref{(ND.routine)} 
should be available.  
This is certainly the case when $p$ happens to be supermodular and $b$ submodular, 
since then $\mu b-p$ is submodular when $\mu \geq 0$
and we can use any submodular function minimization subroutine
(which we abbreviate as {\bf submod-minimizer}).

In several applications, the requested general purpose
submod-minimizer can be superseded by a direct and more efficient
algorithm such as the one for network flows or for matroid partition.
Subroutine  \eqref{(ND.routine)} is also available in the
more general case (needed in applications) 
when the function $p'$ defined by $p'(X):=p(X)-\mu b(X)$ 
is only crossing supermodular
(meaning that the supermodular inequality is expected only 
for pairs $X,Y$ of subsets with $X \cap Y \ne \emptyset$ and $X \cup Y \ne S$).
Indeed, for a given ordered pair of elements $s,t\in S$, 
the restriction of $p'$ to subsets
containing $s$ and avoiding $t$ is fully supermodular,
and therefore we can apply a submod-minimizer 
to each of the $n(n-1)$ ordered pairs $(s,t)$ to get the requested maximum of $p'$.

We call a value $\mu $ \ {\bf good} \ if $\mu b(X)\geq p(X)$ \ 
[i.e., $p(X)-\mu b(X)\leq 0$] \ for every $X\subseteq S$.  
A value that is not good is called {\bf bad}.  
Clearly, if $\mu $ is good, then so is every integer larger than $\mu $. 
We assume that
\begin{equation}
\mbox{$p(X)\leq 0$ \quad whenever \  $b(X)=0$},
\label{(vanjomu)} 
\end{equation}
which is equivalent to requiring that there is a good $\mu$.
We also assume that
\begin{equation}
\mbox{there exists a subset \  $Y\subseteq S$ \ with \  $p(Y)>0$}, 
\label{(0isbad)} 
\end{equation}
which is equivalent to requiring that 
the value $\mu =0$ is bad.
Our goal is to compute the minimum $\mu_{\rm min}$ of the good integers.  
This number is nothing but the maximum of
  \ $\lceil {p(X) / b(X)}\rceil $ \ 
over the subsets of $S$ with $b(X)>0$.

The algorithm starts with the bad $\mu_{0}:=0$.  
Let 
\[
  X_{0} \in \arg\max \{ p(X)-\mu_{0}b(X) :  \ X\subseteq S \},
\]
that is, \ $X_{0}$ is a set maximizing the function $p(X)-\mu_{0}b(X)=p(X)$.  
Note that the badness of $\mu_{0}$ implies that
$p(X_{0}) > 0$.  Since, by the assumption, there is a good $\mu $, 
it follows that $\mu b(X_{0}) \geq p(X_{0}),$ and hence $b(X_{0})>0$.

The procedure determines one by one a series of pairs $(\mu_{j},X_{j})$
for subscripts $j=1,2,\dots $ where each integer $\mu_{j}$ 
is a tentative candidate 
for $\mu $ while $X_{j}$ is a non-empty subset of $S$
with $b(X_{j}) > 0$.
Suppose that the pair $(\mu_{j-1},X_{j-1})$
has already been determined for a subscript $j\geq 1$.  
Let $\mu_{j}$ be the smallest integer 
for which $\mu_{j}b(X_{j-1})\geq p(X_{j-1})$,
that is,
\[
 \mu_{j}:  =\llceil {p(X_{j-1}) \over b(X_{j-1})}\rrceil . 
\]

If $\mu_{j}$ is bad, that is, if there is a set $X\subseteq S$ with
$p(X) -\mu_{j}b(X) > 0$, then let 
\[
  X_{j} \in \arg\max \{ p(X)-\mu_{j}b(X) :  \ X\subseteq S \},
\] 
that is, \ $X_{j}$ is a set
maximizing the function \ $p(X)-\mu_{j}b(X)$.  
(If there are more than one maximizing set, 
we can take any).  
Since $\mu_{j}$ is bad, $X_{j}\not =\emptyset $ and $p(X_{j}) - \mu_{j}b(X_{j})>0$,
which implies $b(X_{j}) > 0$ by the assumption \eqref{(vanjomu)}.

\begin{claim}  \label{betano} 
If $\mu_{j}$ is bad for some subscript $j\geq 0$, then $\mu_{j} < \mu_{j+1}$.  
\end{claim}

\Proof 
The badness of $\mu_{j}$ means that 
$p(X_{j})-\mu_{j}b(X_{j}) > 0$ from which
\[
\mu_{j+1} = \llceil {p(X_{j}) \over b(X_{j})}\rrceil 
  = \llceil {p(X_{j})-\mu_{j}b(X_{j}) \over b(X_{j}) }\rrceil + \mu_{j} 
\ > \  \mu_{j}.  
\]
\vspace{-1.5\baselineskip} \\
\finbox

\medskip

Since there is a good $\mu $ and the sequence $\mu_{j}$ is 
strictly monotone increasing by Claim \ref{betano}, 
there will be a first subscript $h\geq 1$ for which $\mu_{h}$ is good.  
The algorithm terminates by outputting this $\mu_{h}$ 
(and in this case $X_{h}$ is not computed).

\begin{theorem}   \label{Msteps} 
If $h$ is the first subscript during the run of the algorithm 
for which $\mu_{h}$ is good, then 
$\mu_{\rm min}=\mu_{h}$ 
(that is, $\mu_{h}$ is the requested smallest good $\mu $-value)
and $h\leq M$, where $M$ denotes the largest value of $b$.  
\end{theorem}

\Proof 
Since $\mu_{h}$ is good and $\mu_{h}$ is the smallest integer for
which $\mu_{h}b(X_{h-1})\geq p(X_{h-1})$, 
the set $X_{h-1}$ certifies that no good integer $\mu $ 
can exist which is smaller than $\mu_{h}$,
that is, $\mu_{\rm min}=\mu_{h}$.

\begin{claim}  \label{Xno} 
If $\mu_{j}$ is bad for some subscript $j\geq 1$, then $b(X_{j-1}) > b(X_{j})$.  
\end{claim}

\Proof
As $\mu_{j}$ \ ($= \lceil {p(X_{j-1}) / b(X_{j-1})}\rceil$) is bad, 
we obtain that
\begin{align*}
p(X_{j})-\mu_{j}b(X_{j}) >0 & = p(X_{j-1}) - { p(X_{j-1}) \over b(X_{j-1})} b(X_{j-1}) 
\\ &
\geq p(X_{j-1}) - \llceil {p(X_{j-1}) \over b(X_{j-1})}\rrceil b(X_{j-1})
= p(X_{j-1}) - \mu_{j}b(X_{j-1}) ,
\end{align*}
from which we get
\begin{equation} \label{NDprfA}
p(X_{j}) - \mu_{j}b(X_{j}) > p(X_{j-1}) - \mu_{j}b(X_{j-1}).
\end{equation}
\noindent
Since $X_{j-1}$ maximizes $p(X) - \mu_{j-1}b(X)$, we have
\begin{equation} \label{NDprfB}
p(X_{j-1}) - \mu_{j-1}b(X_{j-1}) 
\geq p(X_{j}) - \mu_{j-1}b(X_{j}).  
\end{equation}
\noindent
By adding up 
\eqref{NDprfA} and \eqref{NDprfB},
we obtain
\[
 (\mu_{j} - \mu_{j-1})b(X_{j-1}) > (\mu_{j} - \mu_{j-1})b(X_{j}).
\]
\noindent
As $\mu_{j}$ is bad, so is $\mu_{j-1}$, and hence, 
by applying Claim \ref{betano} to $j-1$ in place of $j$, we obtain that
$\mu_{j} > \mu_{j-1}$, from which we arrive at 
$b(X_{j-1}) > b(X_{j})$,
as required.  
\finbox

\medskip

Claim \ref{Xno} implies that 
$M \geq b(X_{0}) > b(X_{1})> \cdots >b(X_{h-1})$,
from which $1\leq b(X_{h-1}) \leq M -(h-1)$, and hence $h\leq M$ follows.  
This completes the proof of Theorem~\ref{Msteps}.
\finbox \finboxHere 
\medskip

\begin{remark} \rm \label{RMnewtdinksubmax}
The presented variant of the Newton--Dinkelbach algorithm to maximize 
$\left\lceil {p(X) / {b(X)}}\right\rceil$ over subsets $X$ with $b(X)>0$ 
has been shown to be a polynomial algorithm 
for a supermodular function $p$ and a non-negative 
and submodular function $b$
 when the largest value $M$ of $b$ is bounded by a polynomial of $\vert S\vert $,
provided that the seemingly artificial
additional assumptions
in \eqref{(vanjomu)} and \eqref{(0isbad)} hold true.  
However, there is a tiny but sensitive issue here, indicating
that, without these additional
assumptions, the Newton--Dinkelbach (or any other) algorithm cannot 
solve this maximization problem.  
To see this, consider the special case when $b$ is a (finite-valued) submodular function 
which is strictly positive on every non-empty subset, and 
let $N$ be an integer upper bound for 
the squared maximum value of $b$.  
Let $p$ be the function that is identically equal to $-N$ 
except for $p(\emptyset)=0$.
Then $p$ is supermodular.  
Now maximizing $\left\lceil {p(X) / {b(X)}}\right\rceil$ 
is the same as minimizing 
$\lfloor {N / {b(X)}} \rfloor $, 
which is equivalent to maximizing $b(X)$, a well-known NP-hard problem,
even in the case when the maximum value $M$ of $b$ 
is bounded by a polynomial of $|S|$. 
Note that for this special choice of $p$ and $b$, 
the hypothesis \eqref{(0isbad)} fails to hold.  
\finbox
\end{remark}

\subsection{Computing $\beta_{F}$ in strongly polynomial time}
\label{SCminlargeflow}

We describe a strongly polynomial algorithm to compute 
$\beta := \beta_{F}$ in \eqref{(betaF)},
which is the smallest integer for which 
$\odotZ{Q}$ has an element $z$ satisfying $z(e)\leq \beta$ for every edge $e\in F$.
We shall apply the Newton--Dinkelbach algorithm described 
in Section~\ref{SCpbmaxND}
to a supermodular function $p'$  and a submodular function $b$ 
to be defined in \eqref{supmodfnforND} and \eqref{submodfnforND}.

As before, we suppose that there is an $(f,g)$-bounded $m$-flow, and
also that $F$ contains no $(f,g)$-tight edges.  Our first goal is to
find the smallest integer $\beta$ such that by decreasing $g(e)$ to
$\beta$ for each edge $e\in F$ for which $g(e)>\beta$, 
the resulting $g'$ and the unchanged $f$ continue to meet the inequality $f\leq g'$
and the Hoffman-condition \eqref{(Hoffman)}.
The first requirement implies that $\beta$ 
is at least the largest $f$-value on the edges in $F$, which is denoted by $f_{1}$.

Let $g_{1}>g_{2}>\cdots >g_{q}$ denote the distinct $g$-values of the edges
in $F$, and let $L:=\{e\in F:  g(e)=g_{1}\}$.  
Let $\beta_{1}:=\max \{f_{1},g_{2}\}$.

By an $m$-flow feasibility computation, we can check whether the
$g$-value $g_{1}$ on the elements of $L$ can be uniformly decreased to
$\beta_{1}$ without destroying \eqref{(Hoffman)}.  
If this is the case,
then either $\beta_{1}=f_{1}$ in which case a tight edge arises in $F$
and we can remove this tight edge from $F$, or $\beta_{1}=g_{2}$ in which
case the number of distinct 
$g$-values becomes one smaller.
Clearly, as the total number of distinct $g$-values
in $F$ is at most $\vert F\vert $, 
this kind of reduction may occur at most $\vert F\vert $ times.

Therefore, we are at a case when $g_{1}$ cannot be decreased to $\beta_{1}$ 
without violating \eqref{(Hoffman)}.  
Let us try to figure out the
lowest integer value $\beta$ to which $g_{1}$ can be decreased without
violating \eqref{(Hoffman)}.

Recall that $L=\{e\in F:  g(e)=g_{1}\}$ and let $A_{0}:=A-L$
 (that is, $A_{0}$ is the complement of $L$ with respect to the whole edge-set $A$). 
 Let $g'$ denote the function arising from $g$ by reducing
$g(e)$ on the elements of $L$ (where $g(e)=g_{1}$) to $\beta_{1}$.
Since $g'\geq f$ holds and $\varrho_{g'}-\delta_{f}$ is submodular, 
the set-function $p'$ on $V$ defined by
\begin{equation} \label{supmodfnforND}
  p'(Z):  = \widetilde m(Z) - \varrho_{g'}(Z) + \delta_{f}(Z)
\end{equation}
is supermodular.
Define a submodular function $b$ on $V$ by
\begin{equation} \label{submodfnforND}
b(Z):= \varrho_{L}(Z) . 
\end{equation}
Note that the maximum of $b$ is bounded by
a polynomial of the size of the digraph, 
and hence the variant of the Newton--Dinkelbach algorithm
described above is strongly polynomial in this case.

Since $g_{1}$ in the present case cannot be decreased to $\beta_{1}$
without violating \eqref{(Hoffman)}, there is a subset $Z\sp{*}$
violating $\varrho_{g'}(Z) - \delta_{f}(Z) \geq \widetilde m(Z)$, 
or for short, $p'(Z\sp{*})>0$.  

We say that a non-negative integer $\mu $ is {\bf good} 
if it meets the requirement that after increasing uniformly $g(e)=\beta_{1}$ 
by $\mu$ on the edges $e\in L$, Hoffman's condition should hold.  
Our problem to find the smallest $\beta $ is equivalent to computing the
smallest good $\mu $. This is definitely positive since the existence
of $Z\sp{*}$ implies that $\mu =0$ is not good.

\begin{claim}
A positive integer $\mu $ is good if and only if 
\begin{equation}
 \mu b(Z)\geq p'(Z) \quad \hbox{\rm for every $Z\subseteq V$}. 
\label{(mubZ)} 
\end{equation}
\end{claim}

\Proof By definition, $\mu $ is good precisely if 
$$\mu \varrho_{L}(Z) + \varrho _{g'}(Z) - \delta _f(Z) \geq \widetilde m(Z)$$
for every $Z\subseteq V$, which is just equivalent to \eqref{(mubZ)}.
\finbox 
\medskip

The original $g$ meets \eqref{(Hoffman)}, 
meaning that $\varrho_{g}-\delta_{f}\geq \widetilde m$, 
which is equivalent to 
\[
(g_{1}-\beta_{1}) \varrho_{L}(Z) + \varrho_{g'}(Z) -  \delta_{f}(Z) 
   = \varrho_{g}(Z)- \delta_{f}(Z) \geq \widetilde m(Z)
\]
holds for every $Z\subseteq V$.  
This shows that $\mu =g_{1} - \beta_{1}$ is good, 
and our problem requires finding the smallest good $\mu $. 
Since $b$ is submodular, $p'$ is supermodular, and we have 
$\max \{b(Z):Z\subseteq V\} \leq \vert L\vert \leq \vert A\vert $, 
we can apply the 
Newton--Dinkelbach algorithm 
described in Section~\ref{SCpbmaxND}
to this case.

That algorithm needs the subroutine \eqref{(ND.routine)}
to compute a subset of $V$
maximizing $p'(Z) - \mu b(Z)$ \ ($Z\subseteq V$) for any fixed integer
$\mu \geq 0$.  
This subroutine is applied at most $M$ times,
 where $M$ denotes the largest value of $b$.  
Since the largest value of $b$ is at most $\vert A\vert $,
the subroutine 
\eqref{(ND.routine)}
is applied at most $\vert A\vert $ times.  
Furthermore,
by the definition of $p'$ and $b$, the equivalent subroutine to minimize 
\[
 \mu b(Z)-p'(Z)= \mu \varrho_{L}(Z) + \varrho_{g'}(Z)
  -\delta_{f}(Z) - \widetilde m(Z)
\]
can be realized with the help of a straightforward reduction 
to a max-flow min-cut computation in a related 
edge-capacitated digraph on node-set $V\cup \{s,t\}$ with
extra source-node $s$ and sink-node $t$.

Therefore, by relying on an efficient max-flow computation, the
smallest $\mu $ can be computed in strongly polynomial time, and hence
the smallest $\beta \ (= \beta_{1}+ \mu )$ is available for which
$\beta >\beta_{1}=\max \{f_{1},g_{2}\}$ 
and the value $g_{1}$ can be reduced to $\beta$ 
on the edges in $L$ without violating \eqref{(Hoffman)}.

\section{Summary of the algorithm}
\label{SCalgsum}

In this section, we summarize the algorithmic framework discussed in
previous sections.  
We emphasize that each part of the algorithm below is strongly polynomial.  
The input of the algorithm is a digraph $D=(V,A)$, 
integral bounding functions $f\leq g$ on $A$, 
a (finite-valued) integral function $m$ 
with $\widetilde m(V)=0$,
and a subset $F\subseteq A$ of edges, as described in Theorem~\ref{findecmin}.  
Let $Q=Q(f,g;m)$ denote the set of $(f,g)$-bounded $m$-flows, 
while $\odotZ{Q}$ is the set of integral elements of $Q$.

\medskip

{\bf Part~1} of the algorithm decides whether $\odotZ{Q}$ is empty or not.  
This can be done with an adaptation of a max-flow min-cut algorithm.  
So we assume henceforth that $\odotZ{Q}$ is non-empty.

\medskip

If $F=\emptyset $,
then the algorithm terminates with the conclusion 
that every member of $\odotZ{Q}$ is $F$-dec-min.  
So we assume henceforth that $F$ is non-empty.

\medskip

{\bf Part~2} of the algorithm decides whether $\odotZ{Q}$ has an $F$-dec-min element.  
The answer is obviously yes when $f$ and $g$ are finite-valued on $F$.  
In the general case, Part~2 can be realized by the algorithm described 
in Section~\ref{veges},  
which was based on Theorem~\ref{findecmin} 
and Corollary~\ref{fgveges}.  
Part~2 may terminate in two ways.  
In the first one, it outputs a di-circuit $C$ in $D\sp{\infty}$ 
(defined in \eqref{(Ainfdef)}) intersecting $F$.  
Such a di-circuit certifies that no $F$-dec-min element exists.  
In this case, the algorithm terminates with the conclusion 
that $\odotZ{Q}$ has no $F$-dec-min element.  
The other possible output of Part~2 is a new bounding pair $(f',g')$
(described in Corollary~\ref{fgveges}) for which the set 
of $F$-dec-min elements of $\odotZ{Q}(f,g;m)$ is equal to the set 
of $F$-dec-min elements of $\odotZ{Q}(f',g';m)$, 
where $f'$ and $g'$ are finite-valued on $F$.  
In this case there is an $F$-dec-min element of $\odotZ{Q}$.  
Henceforth, we can assume for the remaining parts of the algorithm 
that $f$ and $g$ themselves are finite-valued on $F$.

\medskip 
Suppose that Part~2 is finished.  
In the next parts of the algorithm we need the operation of $F$-reductions.  
\medskip

{\bf $F$-reductions and termination} \ During its run, 
the algorithm carries out edge-tightening steps.  
Such a step (by its definition) does not make necessarily 
an edge $e\in F$ tight, 
but when it does, we carry out an {\bf $F$-reduction} (or an {\bf $F$-reducing step}) 
which is simply the replacement of $F$ by $F-e$.  
An $F$-reduction does not change the set of $F$-dec-min elements. 
If $F$ becomes empty here, 
the whole algorithm terminates with the current bounding pair $(f\sp{*},g\sp{*})$.  
The number of $F$-reductions is at most $|F| \leq |A|$. 
After an application of $F$-reduction, we can assume that the updated 
$F$ contains no tight edges and $F$ is non-empty.

\medskip

{\bf Part~3} of the algorithm computes $\beta := \beta _{F}$ 
defined in \eqref{(betaF)}, 
which is the smallest integer 
for which $\odotZ{Q}$ has an element $z$ satisfying $z(e)\leq \beta$ 
for every edge $e\in F$.  
This is done with the help of the discrete variant of
the Newton--Dinkelbach algorithm in Section~\ref{SCpbmaxND}.  
If we reduce $g(e)$ to $\beta$ for each edge $e\in F$ with $g(e)>\beta $, 
then the set of $F$-dec-min elements does not change.
Therefore we assume henceforth that $\beta =\max\{g(e):e\in F\}$.  
For each edge $e\in F$ with $f(e)=g(e)$, we carry out an $F$-reducing step.  
Let $L:=\{e\in F: g(e)=\beta \}$.  
Note that $L \ne \emptyset$ and $f(e)<g(e)=\beta$ for each $e\in L$, 
and hence the conditions in \eqref{(hypo)} hold.

\medskip

Part~3 finishes by outputting $\beta$. 
Recall that an element $z\in \odotZ{Q}$ 
is said to be pre-dec-min on $F$ if the number $\mu$ of
edges $e\in L$ with $z(e)=g(e)$ ($=\beta$) is minimum.  
Theorem~\ref{Fdecmin}  
states the existence of a certain chain ${\cal C}'$ of subsets of $V$,
called a dual optimal chain,
which provides a certificate 
for an element $z\in \odotZ{Q}$ to be pre-dec-min on $F$.  
In what follows, the algorithm shall apply iteratively Part~4.

\medskip

{\bf Part~4} first computes a dual optimal chain ${\cal C}'$ 
by the algorithm described in Section~\ref{algo5}.  
Next, we consider the updated bounds $(f',g')$ defined 
in \eqref{(f'g'def1)}  and \eqref{(f'g'def2)}  
with reference to ${\cal C}'$.  
As in the proof of Theorem~\ref{reduction}, let
$L'\subseteq L$ consist of those elements of $L$
that enter at least one member of ${\cal C}'$, and let $F':=F-L'$.  
If $F'=\emptyset $, then the whole algorithm terminates 
with the conclusion that the pair $(f\sp{*},g\sp{*})$ defined by $f\sp{*}:=f'$ and
$g\sp{*}:=g'$ meets the requirement of Theorem~\ref{MAIN}.  
If $F'\ne \emptyset $, then we iterate Part~4 for $(f,g):=(f',g')$ and $F:=F'$.
Clearly, the algorithm terminates after at most $|A|$ iterations.

\medskip

If necessary, we can compute a vector-potential certificate, described
by Property (C) in Theorem~\ref{gall.main}, from the pair $(f\sp{*},g\sp{*})$
computed above, as follows.

\medskip

{\bf Part~5}  of the algorithm computes a vector-potential $\underline{\pi}$
for a given $F$-dec-min element $z$ of $\odotZ{Q}$.  
To this end, consider the $k$-dimensional cost-vector 
$\underline{c}$ defined in \eqref{(cdef)}.
By Theorem~\ref{gall.main} and Lemma~\ref{LMekvi},
$\underline{c}$ is conservative, implying that 
there exists a $\underline{c}$-feasible potential-vector $\underline{\pi}$, and this
can actually be computed by Remark~\ref{RMpotveccomp1}.

\section{Remarks on two related problems}
\label{SCconcrem}

\subsection{Fractional dec-min flows}
\label{SCegflfrac}

While we have so far been concerned exclusively with integral flows, it 
is also natural to consider decreasing minimality among real-valued (or fractional) flows
with respect to a specified subset $F$ of edges.
Indeed, the seminal work of Megiddo \cite{Meg74}, \cite{Meg77}
dealt with this continuous (fractional) case when $F$ is the set of edges leaving a source node.
In the following we briefly describe how our structural results 
(Theorems \ref{MAIN},  \ref{gall.main}, and \ref{findecmin})
for the discrete case 
can be adapted to real-valued (fractional) flows.

Let $D=(V,A)$ be a digraph and $F\subseteq A$ 
a non-empty subset of edges.  
Let $m:V \rightarrow {\bf R} $ be a function on $V$ 
with $\widetilde m(V)=0$, 
and let 
$f:A \rightarrow {\bf R}\cup \{-\infty \}$ and
$g:A \rightarrow {\bf R}\cup \{+\infty \}$ be bounding functions on $A$
such that there is an $(f,g)$-bounded $m$-flow in $D$.  
Let $Q=Q(f,g;m)$ denote the set of $(f,g)$-bounded $m$-flows, where
$Q$ is a non-empty subset of ${\bf R}\sp{A}$
consisting of real vectors.
We are interested in $F$-decreasing minimality among members of $Q$.

Concerning the existence of an $F$-dec-min element of $Q$,
we have the following theorem, which is the continuous counterpart of Theorem~\ref{findecmin}.

\begin{theorem}   \label{findecmin-R} 
There exists a (possibly fractional)
$F$-dec-min $(f,g)$-bounded $m$-flow 
if and only if 
there is no di-circuit $C$ with $C\cap F\not =\emptyset$
in the digraph $D\sp{\infty }=(V,A\sp{\infty} )$ defined by \eqref{(Ainfdef)}.
\end{theorem}  
\Proof
The proof is essentially the same as that of Theorem~\ref{findecmin}.
The only difference is that 
the definition of $z'(uv):=z(uv) \pm 1$ in \eqref{existdecminprf1}
should be changed to 
$z'(uv) := z(uv) \pm \delta$
using an arbitrary positive number $\delta > 0$.
\finbox 
\medskip

The characterizations of an $F$-dec-min flow
 for the discrete case
in terms of an improving di-circuit and a potential-vector
(Theorem~\ref{gall.main})
can be adapted to the continuous case as follows.
For a real-valued flow $x: A \rightarrow {\bf R}$ 
we consider the standard auxiliary graph $D_{x}$,
introduced at the beginning of Section~\ref{SCgall}.
The expressions
\eqref{(f'g'def22-R)}, \eqref{(z*def-R)}, and \eqref{(cdef-R)} below
are the continuous counterparts of 
\eqref{(f'g'def22)}, \eqref{(z*def)}, and \eqref{(cdef)}, respectively.

A di-circuit $C$ of $D_{x}$ is called 
{\bf $x$-improving on $F$} (or just $x$-improving) 
if 
there exists a positive number $\delta$
such that $x'$ defined by
\begin{equation} 
x'(uv):= \begin{cases} 
 x(uv)+ \delta  & \hbox{if $uv$ is a forward edge of $C$,}\ 
\cr
 x(uv)- \delta & \hbox{if $vu$ is a backward edge of $C$,}\ 
\cr
x(uv) & \hbox{otherwise}\ 
\end{cases} 
\label{(f'g'def22-R)} 
\end{equation}
for $uv\in A$
is a member of $Q$ and is decreasingly smaller than $x$ on $F$.
Note that the definition of $D_{x}$ implies that $x'$ is indeed in $Q$
for a sufficiently small $\delta > 0$.

The potential-vector 
$\underline{c}$ is defined as follows.
Let $F_{x}$ denote the subset of $A_{x}$ corresponding to $F$,
and let $F_{\bf f}$ and $F_{\bf b}$
be the sets of forward and backward edges in $F_{x}$.
Using the $\delta >0$ above, define a function $x\sp{*}$ on $F_{x}$ by
\begin{equation} 
x\sp{*}(uv):= \begin{cases} 
x(uv) & \hbox{if $uv\in F_{\bf f}$,}\
\cr 
x(vu)- \delta  & \hbox{if $uv\in F_{\bf b}$.}\ 
\end{cases}
\label{(z*def-R)} 
\end{equation}
\noindent
Denoting by  $\gamma_{1} > \gamma_{2}>\cdots >\gamma_{k}$ 
the distinct values of $x\sp{*}$,
we define a $k$-dimensional vector
$\underline{c}(e)$ for every edge $e$ of $D_{x}$ as follows:
\begin{equation} 
\underline{c}(e):= \begin{cases} 
\ \underline{0}_{k} & \hbox{if $e\in A_{x}-F_{x}$,}\ 
\cr
\ \underline{\varepsilon }_{i} & \hbox{if $e\in F_{\bf f}$ and $x\sp{*}(e)=\gamma_{i}$, }\ 
\cr 
\ -\underline{\varepsilon }_{i} & \hbox{if $e\in F_{\bf b}$ and $x\sp{*}(e)=\gamma_{i}$,}\ 
\end{cases} 
\label{(cdef-R)} 
\end{equation}
where $\underline{\varepsilon}_{i}$ is the $k$-dimensional unit-vector
$(0,\dots ,0,1,0, \dots ,0)$ whose $i$-th component is 1. 
Note that the dimension $k$ is bounded by $2|F|$.

With the modified definitions of an improving di-circuit and a potential-vector, 
the following result can be proved 
by modifying the proof of Theorem~\ref{gall.main} in Section~\ref{SCgall}.

\begin{theorem}   \label{gall.main-R} 
For a (possibly fractional) element $x\in Q=Q(f,g;m)$, 
the following properties are equivalent.

\noindent {\rm (A)} 
\ $x$ is decreasingly minimal on $F$.

\noindent {\rm (B)} 
\ There is no $x$-improving di-circuit in the auxiliary digraph $D_{x}$.

\noindent {\rm (B$'$)} 
\ There is no di-circuit $C$ with $\widetilde {\underline{c}}(C)\prec \underline{0}_{k}$
in the auxiliary digraph $D_{x}$.

\noindent {\rm (C)} 
\ There is a potential-vector function 
$\underline{\pi} $ on $V$ which is $\underline{c}$-feasible in $D_{x}$,
that is, $\underline{\pi} (v) - \underline{\pi} (u) \preceq \underline{c}(uv)$ 
for every edge $uv \in A_{x}$.
\finbox
\end{theorem}

In the discrete case
we have given a description of the set of $F$-dec-min integral $m$-flows
in Theorem~\ref{MAIN}
in terms of a pair of bounding functions
$(f\sp{*},g\sp{*})$.
In the continuous case, 
the flow-values of an $F$-dec-min element of $Q$ 
are uniquely determined on $F$ (see Proposition \ref{PRcnvsetdm} below), and 
therefore, the corresponding statement
reads as follows:
\begin{quote}
There exists a pair $(f\sp{*},g\sp{*})$ of bounding functions on $A$
satisfying
$f(e) \leq f\sp{*}(e) = g\sp{*}(e) \leq g(e)$ for $e\in F$
and $f\sp{*}(e)= f(e)$, $g\sp{*}(e)=g(e)$ for $e\in A-F$, 
such that an $(f,g)$-bounded (real-valued) $m$-flow $x$ is $F$-dec-min 
if and only if $x$ is an $(f\sp{*},g\sp{*})$-bounded $m$-flow.
\end{quote}
\noindent
Although the above statement is rather easy to see,
it will be useful when we want to find a cheapest fractional 
feasible $m$-flow that is dec-min on $F$.
It is of course nontrivial to design an algorithm for finding such $(f\sp{*},g\sp{*})$,
which is left for future investigations.

The statement above shows that an $F$-dec-min element of $Q$, when
restricted to $F$, is unique, which is equivalent to saying that the
dec-min element of the projection of $Q$ to $\RR\sp{F}$ is unique.  
This is, actually, a special case of the following
observation concerning general convex sets.

\begin{proposition}   \label{PRcnvsetdm} 
Let $P$ be a convex subset of  ${\bf R}\sp{n}$.
If a dec-min element of $P$ exists,
it is uniquely determined.
\end{proposition}  
\Proof
Suppose, indirectly, that $x$ and $y$ are distinct dec-min elements of $P$.
Let
$\gamma_{1} > \gamma_{2} > \cdots > \gamma_{k}$
denote the distinct values of the components of $x$ and $y$,
and define
$L_{i}(x) := \{ j : x(j) = \gamma_{i}, \  1 \leq j \leq n \}$
and 
$L_{i}(y) := \{ j : y(j) = \gamma_{i}, \  1 \leq j \leq n \}$
for $i=1,2,\dots, k$.
Let $r$ be the smallest index $i$ such that
$L_{i}(x) \not= L_{i}(y)$.
Since
$|L_{r}(x)| = |L_{r}(y)|$
there exist 
$j' \in L_{r}(x) -  L_{r}(y)$ 
and $j'' \in L_{r}(y) -  L_{r}(x)$,
for which 
$x(j') = \gamma_{r} > y(j')$
and 
$y(j'') = \gamma_{r} > x(j'')$.
This implies that 
$(x+y)/2$ is decreasingly smaller than $x$, 
whereas $(x+y)/2$ is in $P$ by the convexity of $P$.
This is a contradiction.
\finbox
\medskip

\subsection{Relation to convex minimization}
\label{SCconvmin}

The dec-min problem 
is often related to minimization of a convex cost function.
For example, if $Q$ is a base-polyhedron,
an element of $Q$ is dec-min in $Q$
if and only if it is a square-sum minimizer of $Q$ \cite{Fuj80,Fuj05book}.
The corresponding statement is also true in its discrete version 
where $Q$ is an M-convex set \cite{FM21partA}.

However, the equivalence between dec-minimality and 
square-sum minimality fails for network flows.
The following example demonstrates that,
both in integral and fractional cases,
an $F$-dec-min flow is not characterized as a feasible flow
with minimum square-sum of flow-values on $F$.

\begin{example} \rm \label{EXdmsqflow}
Consider $D=(V,A)$ with $F \subseteq A$ 
(see Fig.~\ref{FGflowConv}) defined by 
\begin{align*}
& V:= \{ s_{1}, s_{2};  u_{1}, u_{2}, u_{3}, u_{4};  
   v_{1}, v_{2}, v_{3}, v_{4};  t_{1}, t_{2}  \} ,
\\ & 
F :=\{ u_{1}v_{1},  u_{2}v_{2},  u_{3}v_{3},  u_{4}v_{4}   \},
\\ & 
A := \{ s_{1}u_{1},  s_{1}u_{4},  s_{2}u_{2},  s_{2}u_{3}  \} 
\cup F \cup
 \{ v_{1}t_{1},  v_{3}t_{1}, v_{2}t_{2},  v_{4}t_{2}  \} .
\end{align*}
Let $f(e)=0$ and $g(e)=4$ for all $e \in A$, and define $m: V \to \ZZ$ as follows:
\begin{align*}
& m(s_{1}) =  m(s_{2})= -1; \  m(t_{1}) = m(t_{2})= +1 ,
\\ & 
  m(u_{1}) = -2, \ m(u_{2}) = -2, \  m(u_{3}) = -3, \ m(u_{4}) = 0,  
\\ & 
  m(v_{1}) = +2, \ m(v_{2}) = +2, \  m(v_{3}) = +3, \  m(v_{4}) = 0.
\end{align*}

\begin{figure}\begin{center}
\includegraphics[height=35mm]{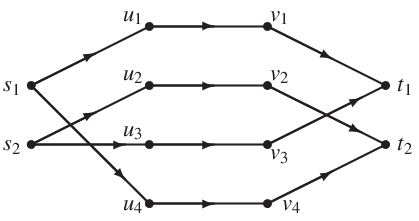}
\caption{Graph $D=(V,A)$ in Example \ref{EXdmsqflow}}
\label{FGflowConv}
\end{center}\end{figure}

There are (precisely) two integral feasible flows, say,
 $x_{1}$ and $x_{2}$,
each corresponding to a pair of disjoint paths from 
$\{ s_{1},s_{2} \}$ to $\{ t_{1},t_{2} \}$,
with additional flows on $F$ required by the condition
$m(u_{i}) =  -m(v_{i})$ for $i=1,2,3,4$.
Their flow-values on $F$ are given by
\[
 x_{1}|_{F} = (2,2,3,0) +  (1,1,0,0) = (3,3,3,0),
\quad 
x_{2}|_{F} = (2,2,3,0) +  (0,0,1,1) = (2,2,4,1),
\]
where $x_{1}$ is the unique $F$-dec-min integral flow.
Nevertheless, $x_{1}$ has a larger square-sum on $F$
than that of $x_{2}$;
the square-sum of $x_{1}|_{F}$ is $27$
and that of $x_{2}|_{F}$ is $25$.

In the fractional (or continuous) case,
the feasible flows are precisely the convex combinations of 
$x_{1}$ and $x_{2}$.
That is,
$x\sp{(\lambda)} = \lambda x_{1} + (1-\lambda) x_{2}$ with $0 \leq \lambda \leq 1$,
and 
\begin{equation}
 x\sp{(\lambda)}|_{F}  
   = \lambda (3,3,3,0) + (1-\lambda) (2,2,4,1) 
   = (2+\lambda, 2+\lambda,4-\lambda,1-\lambda).
\label{fracflowxt}
\end{equation}
This shows that $x_{1} = x\sp{(1)}$ is the unique $F$-dec-min fractional flow.
The square-sum of components of $ x\sp{(\lambda)}|_{F}$ is equal to 
$4 \lambda\sp{2} - 2 \lambda + 25$,
which is minimized at $\lambda=1/4$.
We have
$ x\sp{(1/4)}|_{F} = ( 9/4, 9/4, 15/4, 3/4)$,
which is decreasingly larger than $x\sp{(1)}|_{F}=(3,3,3,0)$. 
Thus, the minimality of square-sum on $F$ 
does not characterize $F$-dec-minimality
even in the fractional case.
\finbox
\end{example}

\medskip

The above example implies, in particular, that an $F$-dec-min fractional flow
cannot be obtained by applying the (strongly polynomial) algorithm
of V{\'e}gh \cite{Veg16} for quadratic-cost fractional flows.

Although the above example denies the use of a quadratic cost function
for the dec-min flow problem, there remains the possibility
of using a more general convex function 
to formulate the dec-min flow problem.
However, the following example indicates that,
in the fractional case, the dec-min flow problem
cannot be formulated as a minimum-cost flow problem 
for any choice of a separable convex objective.

\begin{example} \rm \label{EXdmconvflowR}
Let $\varphi$ be an arbitrary strictly convex (smooth) function on $\RR$.
Referring to the expression \eqref{fracflowxt} of  $x\sp{(\lambda)}|_{F}$, we consider
\[
 \Phi(\lambda) := \varphi(2+\lambda) + \varphi(2+\lambda) 
   + \varphi(4-\lambda) + \varphi(1-\lambda),
\]
which is a separable convex function in the components of $x\sp{(\lambda)}|_{F}$.
Recall that $\lambda=1$ corresponds to $(3,3,3,0)$, which is dec-min
among the vectors $x\sp{(\lambda)}|_{F}$
with $0 \leq \lambda \leq 1$.  
The derivative of $\Phi$ at $\lambda=1$ is positive. Indeed, we have
\begin{align*}
 \Phi'(\lambda) & = 2 \varphi'(2+\lambda) - \varphi'(4-\lambda) - \varphi'(1-\lambda),
\quad 
 \Phi'(1)  =  \varphi'(3) - \varphi'(0) > 0.
\end{align*}
This implies that
the $F$-dec-min flow $x\sp{(1)}$
is not a minimizer of the separable convex function
$\sum_{e \in F} \varphi(x(e))$
over all feasible (fractional) flows $x$.
It is emphasized that such discrepancy 
exists for any choice of $\varphi$.
\finbox
\end{example}

\medskip

The discrepancy of dec-min from convex minimizer
demonstrated above
implies, in particular, that 
algorithms for convex cost flows,
such as those described in the book of 
Ahuja, Magnanti, and Orlin \cite{AMO93},
cannot be used directly for fractional dec-min flow problem.
In this connection, we mention that
the fractional dec-min flow problem
can be solved in (weakly) polynomial time by solving a sequence
of linear programs; see Nace and Orlin \cite{NaOr07}.

In contrast to the fractional case,
the dec-min problem for $Q \subseteq \ZZ\sp{n}$ (in general)
can be formulated as a separable convex function minimization,
as discussed in \cite[Section~3]{FM19partIIarXiv}.
In our integral $F$-dec-min flow problem, we can take, for example,
a real-valued cost function
$\sum_{e \in F} |F|^{x(e)}$
for an integral flow $x$.
Here the function 
$\varphi(k) = |F|^{k}$, defined for all integers $k$,
is increasing and strictly convex in the sense that 
$\varphi(k-1) + \varphi(k+1) > 2 \varphi(k)$ $(k \in \ZZ)$.
Such convex formulation enables us to 
solve the dec-min flow problem 
in (weakly) polynomial time
using the approach of  
Hochbaum and Shanthikumar \cite{HS90}.

\section*{Acknowledgments.}
We thank A. J\"uttner and T. Maehara for illuminating the essence 
of the Newton--Dinkelbach algorithm.  
We are particularly grateful to the anonymous referees
for their thoughtful suggestions,
 which provided a substantial help in forming the final version. 
This research was supported through the program ``Research in Pairs''
by the Mathematisches Forschungsinstitut Oberwolfach in 2019.
The two weeks we could spend at Oberwolfach provided 
an exceptional opportunity to conduct particularly intensive research.
The research was partially supported by the
National Research, Development and Innovation Fund of Hungary
(FK\_18) -- No. NKFI-128673,
and by JSPS KAKENHI Grant Numbers JP26280004, JP20K11697.







\end{document}